\RequirePackage{fix-cm}
\documentclass[smallextended]{svjour3}       
\smartqed  
\usepackage{lipsum}
\usepackage{amsfonts}
\usepackage{graphicx}
\usepackage{epstopdf}
\usepackage{algorithmicx}
\usepackage{enumerate}
\usepackage{amsmath}
\usepackage{amssymb}
\usepackage[utf8]{inputenc}
\usepackage{algorithm}
\usepackage{algpseudocode}
\usepackage{color}
\usepackage{caption}
\usepackage{subcaption}
\usepackage{url}
\usepackage{adjustbox}
\usepackage{multirow}

\newtheorem{assumption}[theorem]{Assumption}
\numberwithin{equation}{section}

\newcommand{\argmax}{\operatornamewithlimits{argmax}}
\newcommand{\argmin}{\operatornamewithlimits{argmin}}

\,
\,

\newcommand{\LCal}{\mathcal{L}}

\newcommand{\dom}{\textnormal{dom}}
\newcommand{\intn}{\textnormal{int}}

\def\proof{\par\noindent{\em Proof. }}
\def\endproof{\hfill $\Box$ \vskip 0.4cm}

\makeatletter
\def\BState{\State\hskip-\ALG@thistlm}
\makeatother

\newcommand{\be}{\begin{equation}}
\newcommand{\ee}{\end{equation}}
\newcommand{\ba}{\begin{array}}
	\newcommand{\ea}{\end{array}}
\newcommand{\bpm}{\begin{pmatrix}}
	\newcommand{\epm}{\end{pmatrix}}

\newcommand{\Rset}{\mathbb{R}}

\newcommand{\ie}{\textit{i}.\textit{e}.}
\newcommand{\eg}{\textit{e}.\textit{g}.}
\ifpdf
\DeclareGraphicsExtensions{.eps,.pdf,.png,.jpg}
\else
\DeclareGraphicsExtensions{.eps}

\usepackage{amsopn}

\fi

\numberwithin{theorem}{section}
\begin{document}

\title{Convergence of the Augmented Decomposition Algorithm
}


\author{Hongsheng Liu         \and
        Shu Lu 
}


\institute{Hongsheng Liu \at
              Department of Statistics and Operations Research, University of North Carolina at Chapel Hill \\
              \email{hsliu@live.unc.edu}           
           \and
           Shu Lu\at{Department of Statistics and Operations Research, University of North Carolina at Chapel Hill\\ \email{shulu@email.unc.edu}}
}

\date{Received: date / Accepted: date}

\maketitle

\begin{abstract}
We study the convergence of the Augmented Decomposition Algorithm (ADA) proposed in \cite{rockafellar2017PSPA} for solving multi-block separable convex minimization problems subject to linear constraints. We show that the global convergence rate of the exact ADA is $o(1/\nu)$ under the assumption that there exists a saddle point.
We consider the inexact Augmented Decomposition Algorithm (iADA)  and  establish global and local convergence results under some mild assumptions, by providing a stability result for the maximal monotone operator $\mathcal{T}$ associated with the perturbation from both primal and dual perspectives. This result implies the local linear convergence of the inexact ADA for many applications such as the \textit{lasso}, total variation reconstruction, exchange problem and many other problems from statistics, machine learning and engineering with $\ell_1$ regularization. 

\end{abstract}

\section{Introduction}\label{sec:Intro}

Consider the following convex optimization problem of minimizing the sum of $K$ separable, potentially nonsmooth convex functions subject to the linear constraints
\begin{equation}\label{P}
\begin{aligned}
& \underset{x}{\text{min}}
& & f(x)=f_{1}(x_1)+\cdots+f_{K}(x_K)\\ 
& \text{ s.t.}
& & Ex=E_1x_1+\cdots+E_Kx_K= q,\\
&
& & x_k\in X_k,\quad k=1,2,\ldots,K, 
\end{aligned}
\end{equation}
where every $f_k$ is a  closed proper convex function (possibly nonsmooth) and each $X_k$ is a closed convex set in $\Rset^{n_k}$. Let $x=(x_1,\ldots,x_K)\in\Rset^n$ be a partition of the variable $x$ and $X=X_1\times \cdots \times X_K\subset \mathbb{R}^{n_1}\times \cdots \times\mathbb{R}^{n_K}=\mathbb{R}^n$ be the domain of $x$. For the linear constraint, $E=(E_1,\ldots,E_K)\in\Rset^{m\times n}$ is a partition of the matrix $E$ consistent with the partition of $x$ and $q\in\Rset^m$ is a column vector. A linear inequality constraint of the form $Ex\le q$ can be easily transformed to the equality case by introducing a slack variable $x_{K+1}\ge 0$. 

Optimization problems in the form of \eqref{P} arise in many application areas such as signal processing, statistics and machine learning. 
\cite{ma2016alternating} summarizes a list of applications arising from many areas when more than two blocks are involved ($K\ge 3$).

Many decomposition algorithms have been proposed to solve the above optimization problem; see \cite{chan2014distributed,chatzipanagiotis2015augmented,chen1994proximal,eckstein1992douglas,mulvey1992diagonal,spingarn1985applications,tseng1991applications} and references therein. Among them, the ADMM method is perhaps the most popular approach to solve the decomposition problem due to its suitable parallel implementation and outstanding computational performance. When $K=2$, the convergence of the ADMM was well studied in the framework of Douglas-Rachford splitting method \cite{eckstein1992douglas}. The paper \cite{deng2016global} proved the linear convergence of the ADMM when at least one of $f_i(\cdot)$ is strongly convex and $E$ satisfies some additional assumptions. For the $K\ge 3$ case, it was shown in \cite{han2012note} that the global convergence is guaranteed if all objective functions $f_k$ are strongly convex. However, for general convex objective functions, it is acknowledged that the direct extension of the original ADMM may diverge \cite{chen2016direct}. Therefore, most recent researches have been focused on either analyzing problems with additional assumptions or showing the convergence results for variants of the ADMM; see \cite{hong2017linear,wang2013solving}. 

As an alternative to the ADMM algorithm for multi-block convex optimization problems, a new primal-dual algorithm called the \textit{augmented decomposition algorithm} (ADA) was introduced in \cite{rockafellar2017PSPA}. This method is closely related to the decomposition algorithm based on the partial inverses proposed in \cite{spingarn1985applications} but is derived from the proximal saddle point algorithm (PSPA) which is associated with a special primal-dual saddle function. It was shown in \cite{rockafellar2017PSPA} that the algorithm is guaranteed to converge on the basis of convergence results of the proximal point algorithm (PPA) in \cite{rockafellar1976augmented}. What is more exciting is that the calculation of each iteration in PSPA can be carried out in parallel and its parallel implementation leads to the ADA. 

Although the global convergence result for the ADA has been well studied under a general condition, the convergence rate result remained unknown. In the first part of this paper, we focus on the convergence analysis of the ADA applied to problem \eqref{P}. 
For that, we first provide a detailed proof for its convergence. Then, we show the $O(1/\nu)$ convergence rate in an ergodic sense. Finally, we improve the convergence result from $O(1/\nu)$ to $o(1/\nu)$ in a non-ergodic sense. These ideas are inspired by recent works on the ADMM and variants of the proximal method of multiplier \cite{deng2017parallel,he20121,shefi2014rate}.

Then, we consider the inexact ADA (iADA) in the second part. We first establish the global convergence result under certain approximation criteria. Then, under some mild assumptions on the function $f_k$ and the structure of feasible set $X_k$, we show the local linear convergence of the iADA. This work is invoked by recent convergence rate results for the ADMM algorithm in \cite{deng2016global,hong2017linear}. 
However, our proof is different from them in which we show the stability of a maximal monotone operator associated with the saddle function for a variant of \eqref{P}. Denote the Lagrangian function by $L$ for \eqref{P}:
\begin{equation}\label{L:P}
L(x,y)=\begin{cases}
f(x)+\langle Ex-q,y\rangle,&\quad \forall (x,y)\in X\times\Rset^m,\\
\infty,&\quad \forall x\notin X.
\end{cases}
\end{equation}
The corresponding maximal monotone operator $\mathcal{T}_L$  \cite{rockafellar1976augmented} is defined by 
\begin{equation}\label{T_L}
\mathcal{T}_L(x,y)=\{(u,v)|(u,-v)\in\partial L(x,y)\}
\end{equation}
where $\partial L(x,y)$  denotes the subgradient of the convex-concave function $L$.
The inverse of $\mathcal{T}_L$ is given by
\begin{equation}\label{T_L^{-1}}
\mathcal{T}_L^{-1}(u,v) = \{(x,y)|(u,-v)\in\partial L(x,y)\}.
\end{equation}
 A solution to $(0,0)\in\mathcal{T}_L(x,y)$ is a saddle point of $L$. 
Classical convergence rate results for PPA \cite{rockafellar1976monotone} rely on the assumption that $\mathcal{T}_L^{-1}$ is Lipschitz continuous at $(0,0)$. This result was extended in \cite{luque1984asymptotic} for situations in which  $\mathcal{T}_L^{-1}(0,0)$ is not a singleton and the following holds:
\begin{equation}\label{Lip:T_L}\small
\exists a>0,\quad \exists \delta>0: \quad \forall w\in\mathcal{B}((0,0),\delta), \quad\forall z\in \mathcal{T}_L^{-1}w, \quad dist(z,\mathcal{T}_L^{-1}(0,0))\le a||w||.
\end{equation}
It has been pointed out in many works that understanding the Lipschitzian behavior of $\mathcal{T}_L^{-1}$ at the origin is crucial to the study of the local convergence results for algorithms in the PPA framework; see \cite{cui2017r,han2015linear,li2016highly,liu2012implementable}. For instance, \cite{li2016highly} showed the \emph{metric subregularity} defined in \cite{dontchev2009implicit} of $\mathcal{T}_L$ which is closely related to \eqref{Lip:T_L} under the so-called second order sufficient condition. However, this result inherently requires the solution uniqueness for problem \eqref{P}.  
Compared with those assumptions, our assumptions in this part mainly rely on the polyhedral property of the feasible set $X$ and the optimal solution set for \eqref{P} needs not to be a singleton. Our proof is based on Robinson's celebrated work on the error bound result for polyhedral multifunctions \cite{robinson1981some} and uses some ideas in the analysis for the satisfaction of a certain error bound condition in \cite{hong2017linear,luo1993convergence}. 

\textbf{Organization }The remainder of this paper is organized as follows. Section \ref{sec:ADA} first summarizes the basic idea of the proximal saddle point algorithm and its implementation, the ADA. Then, we show the convergence result for the ADA and compare it with the ADMM. In Section \ref{sec:iADA}, we introduce the iADA and make some basic assumptions on the problem \eqref{P} for further discussion. Section \ref{sec:stability} studies the stability results of the maximal monotone operator $\mathcal{T}_{\bar{L}}$. Section \ref{sec:Conv_iADA} establishes the global convergence and local linear convergence rate results of the iADA. Finally, some numerical examples are presented in Section \ref{sec:numerical} to demonstrate the performance the ADA and iADA.

\textbf{Notation } We use $\langle \cdot,\cdot\rangle$ and $||\cdot||$ to denote the standard inner product  and $\LCal_2$-norm in the Euclidean space respectively. For any positive definite matrix $G\in S_{++}^n$ and $x,y\in\Rset^n$, the inner product $\langle x,y\rangle_G $ is defined by $x^TGy$ and its induced norm is denoted by $||\cdot||_G$. For $1\le q\le\infty$, $||\cdot||_q$ represents the $\LCal_q$-norm. For any $E\in\Rset^{m\times n}$, $||E||$ denotes the spectral norm,
\ie, the largest singular value of $E$. For any function $f$, let \dom$f$ be the effective domain of the function $f$ and \intn(\dom$f$) be the interior of \dom$f$.  For any point $x\in\Rset^n$ and a closed convex set $C\subset \Rset^n$, $dist(x,C)=\min_{y\in C}||y-x||$.

\section{Global convergence of the ADA}\label{sec:ADA}
In this paper, we make the following standard assumption. 
\begin{assumption}\label{assu:saddle-point}
	The global minimum of \eqref{P} is attainable and 
	\vspace{-1.5ex}
	\begin{equation}\label{assu:int_1}
	\vspace{-1.5ex}
	\intn(X)\cap \dom f\cap\{x|Ex=q\}\neq \emptyset.
	\end{equation} 
	If $X$ is  polyhedral, an alternative assumption for \eqref{assu:int_1} can be that 
		\vspace{-1.5ex}
	\begin{equation}\label{assu:int_2}
		\vspace{-1.5ex}
	X\cap \intn(\dom f)\cap\{x|Ex=q\}\neq\emptyset.
	\end{equation}
\end{assumption}

Assumption \ref{assu:saddle-point} guarantees the existence of a saddle point of $L$. Namely, there exist $\bar{x}$ and $\bar{y}$ such that
	\vspace{-1ex}
\begin{equation}
	\vspace{-1ex}
\bar{x}\in\argmin_{x\in X}L(x,\bar{y}),\qquad \bar{y}\in\argmax_{y\in\Rset^{m}}L(\bar{x},y).
\end{equation}
The dual function for problem \eqref{P} is
	\vspace{-1ex}
\begin{equation}\label{dual_P}
	\vspace{-1ex}
d(y)=\min_{x\in X}L(x,y)=\min_{x\in X}\{ f(x)+\langle y, Ex-q\rangle\}
\end{equation}
and its associated dual problem is given by
	\vspace{-1.5ex}
\begin{equation}\label{D}
	\vspace{-1.5ex}
\max_{y\in\Rset^m} d(y).
\end{equation}
Let $X^*$ and $Y^*$ be the optimal solution sets of \eqref{P} and \eqref{D} respectively. The set of saddle points for the Lagrangian \eqref{L:P} is given by $X^*\times Y^*$.

\subsection{Augmented Decomposition Algorithm}
Here, we first summarize the basic idea of PSPA and its parallel implementation, the ADA. For that, the original problem \eqref{P} is equivalently transformed into 
\begin{equation}\label{P:w}
	\vspace{-1.5ex}
\begin{aligned}
& \underset{x,w}{\text{min}}
& & f(x)=f_{1}(x_1)+\cdots+f_{K}(x_K)\\
& \text{ s.t.}
& & E_jx_j-w_j = 0 ,\quad j = 1, \ldots, K-1,\\
&
& & E_Kx_K-q-w_K = 0,\\
&  
& &w_1+\cdots+w_K = 0,\\
&
& & x_k\in X_k,\quad k=1,2,\ldots,K.
\end{aligned}
\end{equation}
If $x=(x_1,\ldots,x_K)\in\Rset^n$ is an optimal solution of \eqref{P}, then  $(x,w)=(x_1,\ldots,x_K,E_1x_1,\ldots,E_{K-1}x_{K-1},E_Kx_K-q)$ will be an optimal solution to \eqref{P:w}. Instead of adding a multiplier vector for $w_1+\cdots+w_K = 0$, \cite{rockafellar2017PSPA} introduced $W$ as a subspace of $(\Rset^m)^K$ which is defined as
\begin{equation}\label{W}
W=\{w=(w_1,\ldots,w_K)|w_1+\cdots+w_K=0\}\subset (\mathbb{R}^m)^K.
\end{equation}
The orthogonal complement subspace of $W$ is given by
\begin{equation}\label{W_perp}
W^\perp=\{w=(w_1,\ldots,w_K)|w_1=\cdots=w_K\}\subset (\mathbb{R}^m)^K.
\end{equation}
For any $w=(w_1,\ldots,w_K)\in(\Rset^m)^K$, we use $P_{W^\perp}(w)$ to denote the projection of $w$ onto the subspace $W^\perp$.
In \cite{rockafellar2017PSPA}, the author proposed to add increments $u_{i}\in\Rset^m,i=1,\ldots,K$ to the first $K$ linear constraints in \eqref{P:w} and in addition, add to $w\in W$ a perturbation $v\in W^\perp$. The Lagrangian function associated with this perturbation finally works out in terms of the subspace 
\begin{equation}\label{S}
S =\{(\eta,\zeta)|P_{W^\perp}(\eta)=\zeta\}\subseteq(\mathbb{R}^m)^K\times W^\perp,
\end{equation}
and the functions
\begin{equation}\label{L:sub_v}
L_j(x_j,\eta_j)=\begin{cases}
f_{j}(x_j)+\eta_j\cdot E_jx_j, \text{ if } j=1,\ldots,K-1,\\
f_{K}(x_K)+\eta_K\cdot (E_Kx_K-q), \text{ o.w. }
\end{cases}
\end{equation}
to mean that
\begin{equation}\label{L:final_v}
\bar{L}(w,x,\eta,\zeta)=\begin{cases}
\sum_{j=1}^{K}[L_j(x_j,\eta_j)-\eta_j\cdot w_j], \text{ if }  (w,x)\in W\times X,  (\eta,\zeta)\in S,\\
-\infty, \text{ if } (w,x)\in W\times X, (\eta,\zeta)\notin S,\\
+\infty, \text{ if } (w,x)\notin W\times X.
\end{cases}
\end{equation}
The next lemma shows the relationship between $L(x,y)$ and $\bar{L}(w,x,\eta,\zeta)$.
\begin{lemma}\label{L:relation}
	If $(\bar{w},\bar{x},\bar{\eta},\bar{\zeta})$ is a saddle point of the Lagrangian function  in \eqref{L:final_v}, then $\bar{\eta}_1=\bar{\eta}_2=\cdot=\bar{\eta}_K$ and $(\bar{x},\bar{\eta}_1)$ is a saddle point of \eqref{L:P}. Conversely, let $(\bar{x},\bar{y})$ be a saddle point of \eqref{L:P}, and define $\bar{w}=(E_1\bar{x}_1,\ldots,E_{K-1}\bar{x}_{K-1},E_K\bar{x}_K-q)\in(\Rset^m)^K$, $\bar{\eta}=(\bar{y},\ldots,\bar{y})\in(\Rset^m)^K$ and $\bar{\zeta}=\bar{\eta}$. Then $(\bar{w},\bar{x},\bar{\eta},\bar{\zeta})$ is a saddle point of \eqref{L:final_v}.
\end{lemma}
\proof The dual problem associated with \eqref{L:final_v} is 
\begin{equation}\label{D:final_v}
\max_{(\eta,\zeta)\in S} \{\bar{g}(\eta,\zeta)=\inf_{(w,x)\in W\times X}\bar{L}(w,x,\eta,\zeta)\}
\end{equation}
with its feasible set given by \[\{(\eta,\zeta)|\bar{g}(\eta,\zeta)>-\infty\}\subset S.\] As $w\cdot\eta$ cannot be $\infty$, this implies $\eta_1=\eta_2=\cdot=\eta_K$. As a consequence, the dual problem reduces to 
\begin{equation}
\max_{(\eta,\zeta)\in S} \{\bar{g}(\eta,\zeta)=\inf_{x\in X}f(x)+\langle \eta_1,Ex-q\rangle\} 
\end{equation}
which is equivalent to the dual problem corresponding to \eqref{L:P}. So we can conclude the first part. The second part is similarly based on the above observation for the dual whose proof is omitted here.\endproof

Based on \cite{rockafellar1976augmented}, the proximal method of multipliers is derived by adding both primal and dual proximal terms into the Lagrangian \eqref{L:final_v}. More explicitly, the proximal saddle point algorithm in \cite{rockafellar2017PSPA} can be described as the following:


Generate a sequence of elements $(w^\nu,x^\nu)\in W\times X$ and $(\eta^\nu,\zeta^\nu)\in S$ by letting 
\begin{equation}\label{L:prox_v}\scriptsize
\bar{L}^\nu(w,x,\eta,\zeta)=\bar{L}(w,x,\eta,\zeta)+\frac{\rho}{2}||w-w^\nu||^2+\frac{1}{2c}||x-x^\nu||^2-\frac{1}{2\rho}||\eta-\eta^\nu||^2-\frac{1}{2\rho}||\zeta-\zeta^\nu||^2
\end{equation}
and calculating 
\[
(w^{\nu+1},x^{\nu+1},\eta^{\nu+1},\zeta^{\nu+1})= \text{ unique saddle point of } \bar{L}^\nu(w,x,\eta,\zeta)
\]
with respect to minimizing over $(w,x)\in W\times X$ and maximizing over $(\eta,\zeta)\in S.$
According to \cite{rockafellar1976augmented}, the sequence $(w^{\nu},x^{\nu},\eta^{\nu},\zeta^{\nu})$ generated by the above algorithm from any initial $(w^1,x^1)\in W\times X$ and $(\eta^1,\zeta^1)\in S$ is certain to converge to some saddle point $(\bar{w},\bar{x},\bar{\eta},\bar{\zeta})$ of the Lagrangian $\bar{L}$. With the special structure of the saddle point problem, the calculation of the saddle point in \eqref{L:prox_v} can be carried out in the following parallel algorithm ADA. For simplicity, we denote
\begin{equation}\label{phi_k}\small
\phi_{k,\rho,c}^\nu(x_k)=\begin{cases}
f_k(x_k)+\frac{\rho}{4}||E_kx_k-w_k^\nu+\frac{2}{\rho}y_k^\nu||_2^2+\frac{1}{2c}||x_k-x_k^\nu||_2^2,\quad k=1,\ldots,K-1,\\
f_K(x_K)+\frac{\rho}{4}||E_Kx_K-q-w_K^\nu+\frac{2}{\rho}y_K^\nu||_2^2+\frac{1}{2c}||x_K-x_K^\nu||_2^2,\quad k=K.
\end{cases}
\end{equation}

\begin{algorithm}
	\caption{Augmented decomposition algorithm}\label{exactADA}
	\begin{algorithmic}[1]
		\State Given $w^0\in W, x^0\in X, y^0\in (\Rset^m)^K$
		\For{$\nu=0,1, \dots $}
		\State $x_k^{\nu+1} = \argmin_{x_k\in X_k}\phi_{k,\rho,c}^\nu(x_k), k=1,\ldots,K$
		\State $\eta_{k}^{\nu+1}=\begin{cases}
		y_{k}^\nu+\frac{\rho}{2}[E_kx_k^{\nu+1}-w_{k}^\nu],\text{ if }k=1,\ldots,K-1\\
		y_{K}^\nu+\frac{\rho}{2}[E_Kx_K^{\nu+1}-q-w_{K}^\nu],\text{ if }k=K
		\end{cases}$
				\For  {$k=1,\ldots,K$}
		\State $\zeta_k^{\nu+1}= \frac{1}{K}\sum_{j=1}^{K}\eta_{j}^{\nu+1}$\\
			\State$w_{k}^{\nu+1}=w_{k}^\nu+\frac{1}{\rho}[\eta_{k}^{\nu+1}-\zeta_k^{\nu+1}]$\\
			\State$y_{k}^{\nu+1}=\frac{1}{2}[\eta_{k}^{\nu+1}+\zeta_k^{\nu+1}]$
		\EndFor

		\EndFor
	\end{algorithmic}
\end{algorithm}


\subsection{Convergence of the ADA}
In this subsection, we assume $q=0$ for notational simplicity which will not influence the proofs below. Define  the matrix 
\begin{equation}\label{def:G}
G:=\begin{pmatrix} \rho I_{mK} & & & \\ &\frac{1}{c}I_n& &\\& & \frac{1}{\rho}I_{mK}&\\ & & & \frac{1}{\rho}I_{mK}\end{pmatrix}.
\end{equation}
Hence $G\succ 0$ and $||\cdot||_G$ defines a norm. Let $\hat{u}=(\hat{w},\hat{x},\hat{\eta},\hat{\zeta})$ and $u^\nu=(w^\nu,x^\nu,\eta^\nu,\zeta^\nu)$ where $\hat{u}$ is a saddle point of the Lagrangian function \eqref{L:final_v} and $u^\nu$ is the current iteration point.
The convergence result for ADA was established in \cite{rockafellar2017PSPA} on the basis of convergence results for the classic PPA. Here, we import the result and provide an alternative proof for it.  
\vspace{-3ex}
\begin{theorem}\label{thm:global_conv_ADA}
	Under Assumption \ref{assu:saddle-point}, for any $\rho>0$ and $c>0$, the sequence $\{(w^\nu,x^\nu,y^\nu)\}_{\nu=1}^{\infty}$ generated in $W\times X\times (\Rset^m)^K$ by the ADA from any starting point converges to some $(\bar{w},\bar{x},\bar{y})$ such that
	\begin{enumerate}[(a)]
		\vspace{-3ex}
		\item $(\bar{w},\bar{x})$ solves \eqref{P:w}, hence $\bar{x}$ solves \eqref{P},
		\item $\bar{y}_1=\cdots=\bar{y}_q\in \mathbb{R}^m$, and this common multiplier vector solves \eqref{D}. 
	\end{enumerate}
\end{theorem}
\proof From Assumption \ref{assu:saddle-point} and Lemma \ref{L:relation}, there exists a saddle point $(\hat{w},\hat{x},\hat{\eta},\hat{\zeta})\in W\times X\times S$ of the Lagrangian function \eqref{L:final_v}. For each iteration $\nu+1$, due to the minimax operation on \eqref{L:prox_v}, from the primal perspective, we have the following inequality
\begin{equation}\label{minimax_L_prox_v}
\begin{aligned}
&\sum_{k=1}^{K}f_k(x_k)+\sum_{k=1}^{K}\langle\eta_k^{\nu+1},E_kx_k-w_k\rangle\\
\ge&\sum_{k=1}^{K}f_k(x_k^{\nu+1})+\sum_{k=1}^{K}\langle\eta_k^{\nu+1},E_kx_k^{\nu+1}-w_k^{\nu+1}\rangle +\frac{1}{c}\sum_{k=1}^{K}\langle x_k-x_k^{\nu+1},x_k^\nu-x_k^{\nu+1}\rangle\\
&+\rho\sum_{k=1}^{K}\langle w_k-w_k^{\nu+1},w_k^\nu-w_k^{\nu+1}\rangle
\end{aligned}
\end{equation}
for any $x\in X$ and $w\in W$. Applying $(w,x)=(\hat{w},\hat{x})$ to \eqref{minimax_L_prox_v} and noticing that $E\hat{x}_k=\hat{w}_k,k=1,\ldots,K$, we obtain
\begin{equation}\label{minimax_opt_L_prox_v}
\begin{aligned}
&\min P := \sum_{k=1}^{K}f_k(\hat{x}_k)
\ge\sum_{k=1}^{K}f_k(x_k^{\nu+1})+\sum_{k=1}^{K}\langle\eta_k^{\nu+1},E_kx_k^{\nu+1}-w_k^{\nu+1}\rangle\\
& -\frac{1}{c}\sum_{k=1}^{K}\langle x_k^{\nu+1}-\hat{x}_k,x_k^\nu-x_k^{\nu+1}\rangle-\rho\sum_{k=1}^{K}\langle w_k^{\nu+1}-\hat{w}_k,w_k^\nu-w_k^{\nu+1}\rangle.
\end{aligned}
\end{equation}
Similarly, from the dual perspective and the saddle-point property of $(\hat{w},\hat{x},\hat{\eta},\hat{\zeta})$, the following inequality 
\begin{equation}\label{minimax_opt_L_prox_v_2}
\begin{aligned}
\min P = &\sum_{k=1}^{K}f_k(\hat{x}_k)\le \sum_{k=1}^{K}f_k(x_k^{\nu+1})+\sum_{k=1}^{K}\langle\hat{\eta}_k,E_kx_k^{\nu+1}-w_k^{\nu+1}\rangle\\
\le&\sum_{k=1}^{K}f_k(x_k^{\nu+1})+\sum_{k=1}^{K}\langle\eta_k^{\nu+1},E_kx_k^{\nu+1}-w_k^{\nu+1}\rangle\\
&+\frac{1}{\rho}\sum_{k=1}^{K}\langle \eta_k^{\nu+1}-\hat{\eta}_k,\eta_k^{\nu}-\eta_k^{\nu+1}\rangle+\frac{1}{\rho}\sum_{k=1}^{K}\langle \zeta_k^{\nu+1}-\hat{\zeta}_k,\zeta_k^{\nu}-\zeta_k^{\nu+1}\rangle
\end{aligned}
\end{equation}
holds. Combining the above two inequalities with the following identity
\begin{equation}\label{3identity}
2\langle a-b,c-a\rangle = ||c-b||_2^2-||c-a||_2^2-||b-a||_2^2, 
\end{equation}
we have
\begin{equation}
\begin{aligned}
&\sum_{k=1}^{K}(\frac{1}{c}||x_k^\nu-\hat{x}_k||_2^2+\rho||w_k^\nu-\hat{w}_k||_2^2+\frac{1}{\rho}||\eta_k^\nu-\hat{\eta}_k||_2^2+\frac{1}{\rho}||\zeta_k^\nu-\hat{\zeta}_k||_2^2)\\
-&\sum_{k=1}^{K}(\frac{1}{c}||x_k^{\nu+1}-\hat{x}_k||_2^2+\rho||w_k^{\nu+1}-\hat{w}_k||_2^2+\frac{1}{\rho}||\eta_k^{\nu+1}-\hat{\eta}_k||_2+\frac{1}{\rho}||\zeta_k^{\nu+1}-\hat{\zeta}_k||_2^2)\\
\ge&\sum_{k=1}^{K}(\frac{1}{c}||x_k^{\nu+1}-x_k^\nu||_2^2+\rho||w_k^{\nu+1}-w_k^\nu||_2^2+\frac{1}{\rho}||\eta_k^{\nu+1}-\eta_k^\nu||_2^2+\frac{1}{\rho}||\zeta_k^{\nu+1}-\zeta_k^\nu||_2^2)
\end{aligned}
\end{equation}
which is equivalent with 
\begin{equation}
||u^\nu-\hat{u}||_G^2-||u^{\nu+1}-\hat{u}||_G^2\ge||u^\nu-u^{\nu+1}||_G^2.
\end{equation}
From this inequality, we can easily conclude that
\begin{itemize}
	\item[(i)] $\sum_{\nu=0}^{\infty}||u^\nu-u^{\nu+1}||_G^2<\infty$;
	\item[(ii)] $\{u^\nu=(w^{\nu},x^{\nu},\eta^{\nu},\zeta^{\nu})\}$ lies in a compact region;
	\item[(iii)] $||u^\nu-\hat{u}||_G$ is a monotonically non-increasing sequence and thus converges.	
\end{itemize}

From (ii), by passing to a subsequence if necessary, there exists at least one limiting point of $\{(w^{\nu},x^{\nu},\eta^{\nu},\zeta^{\nu})\}$, denoted as $\bar{u}=(\bar{w},\bar{x},\bar{\eta},\bar{\zeta})$. It follows from (i) that $x^{\nu}-x^{\nu+1}\rightarrow 0$, $w^{\nu}-w^{\nu+1}\rightarrow 0$ and $\eta^{\nu}-\eta^{\nu+1}\rightarrow 0$. The update rule for $w$ implies that $\bar{\eta}_1=\dots=\bar{\eta}_K$ and thus $\bar{y}_1=\dots=\bar{y}_K=\bar{\eta}_1$. Since $\eta_{k}^{\nu+1}= y_{k}^\nu+\frac{\rho}{2}[E_kx_k^{\nu+1}-w_{k}^\nu]$,  $E_k\bar{x}_k=\bar{w}_k$ holds and thus $E\bar{x}=0$ which implies the feasibility of $\bar{x}$. Due to the optimality condition for each block in iteration $\nu+1$, we have
$$
0\in\partial f_k(x_k^{\nu+1})+E_k^T\eta_k^{\nu+1}+\frac{1}{c}(x_k^{\nu+1}-x_k^{\nu})+N_{X_k}(x_k^{\nu+1}), \quad k=1,\ldots,K.
$$
By passing to  the limit, we obtain
$$
0\in\partial f(\bar{x})+E^T\bar{\eta}_1+N_{X}(\bar{x}).
$$ 
As a result, $(\bar{w},\bar{x},\bar{\eta},\bar{\zeta})$ is a saddle point of the Lagrangian function \eqref{L:final_v}. Next, we show the uniqueness of the limit point to complete the proof. Let $\bar{u}^1=(\bar{w}^1,\bar{x}^1,\bar{\eta}^1,\bar{\zeta}^1)$ and $\bar{u}^2=(\bar{w}^2,\bar{x}^2,\bar{\eta}^2,\bar{\zeta}^2)$ be any two different limit points of $u^\nu=(w^{\nu},x^{\nu},\eta^{\nu},\zeta^{\nu})$. By the previous argument, both of them are saddle points of \eqref{L:final_v}. From (iii), we know the existence of the following limits
$$
\lim_{\nu\rightarrow \infty}||u^\nu-\bar{u}^i||_G=\beta_i,\quad i=1,2.
$$
With the following equality
$$
||u^\nu-\bar{u}^1||_G^2-||u^\nu-\bar{u}^2||_G^2=-2\langle u^\nu,\bar{u}^1-\bar{u}^2\rangle_G+||\bar{u}^1||_G^2-||\bar{u}^2||_G^2
$$
and by passing to the limit, we have
$$
\beta_1^2-\beta_2^2=-2\langle \bar{u}^1,\bar{u}^1-\bar{u}^2\rangle_G+||\bar{u}^1||_G^2-||\bar{u}^2||_G^2=-||\bar{u}^1-\bar{u}^2||_G^2
$$
and
$$
\beta_1^2-\beta_2^2=-2\langle \bar{u}^2,\bar{u}^1-\bar{u}^2\rangle_G+||\bar{u}^1||_G^2-||\bar{u}^2||_G^2=||\bar{u}^1-\bar{u}^2||_G^2.
$$
Thus we obtain $||\bar{u}^1-\bar{u}^2||_G=0$ which implies that the sequence $(w^{\nu},x^{\nu},\eta^{\nu},\zeta^{\nu})$ converges to some saddle point of the Lagrangian function \eqref{L:final_v} and hence (a) and (b) hold.\endproof

\subsection{Rate of Convergence} In this subsection, we study the global convergence rate for the ADA. We first show the sublinear convergence result of the ADA in an ergodic sense. The proof follows the same idea as that in \cite{shefi2014rate}.
\begin{theorem}\label{thm:sublinear}
	Let $\{u^\nu=(w^\nu,x^\nu,\eta^\nu,\zeta^\nu)\}$ in $W\times X\times S$ be the infinite sequence generated by the ADA. For any integer $N>0$, define $\tilde{x}_{N}$ by 
		\vspace{-1ex}
	\be
	\vspace{-1ex}
	\begin{aligned}
		\tilde{x}_{N}=\frac{1}{N}\sum_{\nu=1}^{N}x^\nu.
	\end{aligned}
	\ee
	Then for any saddle point $\hat{u}=(\hat{w},\hat{x},\hat{\eta},\hat{\zeta})\in W\times X\times S$ of \eqref{L:final_v},
	\be\label{conv:sublinear}
	\begin{aligned}
		f(\tilde{x}_N)+\langle \hat{\eta}_1,E\tilde{x}_N\rangle - \min P \le  \frac{||\hat{u}-u^0||_G^2}{N}.
	\end{aligned}
	\ee
\end{theorem} 
\proof For any saddle point $(\hat{w},\hat{x},\hat{\eta},\hat{\zeta})\in W\times X\times S$ of the Lagrangian function \eqref{L:final_v}, it follows from \eqref{minimax_opt_L_prox_v} and \eqref{minimax_opt_L_prox_v_2} that 
\begin{equation}\label{decrease:primal-dual}
\begin{aligned}
&||u^\nu-\hat{u}||_G^2-||u^{\nu+1}-\hat{u}||_G^2 \\\ge&||u^\nu-u^{\nu+1}||_G^2+\sum_{k=1}^{K}f_k(x_k^{\nu+1})+\sum_{k=1}^{K}\langle\hat{\eta}_k,E_kx_k^{\nu+1}\rangle-\min P\\\ge&\sum_{k=1}^{K}f_k(x_k^{\nu+1})+\sum_{k=1}^{K}\langle\hat{\eta}_k,E_kx_k^{\nu+1}\rangle-\min P.
\end{aligned}
\end{equation}
Summing \eqref{decrease:primal-dual} for $\nu=0,1,\ldots,N-1$, we obtain
\begin{equation}\label{decrease:primal-dual-sum}
\begin{aligned}
&||u^0-\hat{u}||_G^2\\\ge&\sum_{\nu=0}^{N-1}\{\sum_{k=1}^{K}f_k(x_k^{\nu+1})+\sum_{k=1}^{K}\langle\hat{\eta}_k,E_kx_k^{\nu+1}\rangle\}-N\min P\\\ge& N[f(\tilde{x}_N)+\langle \hat{\eta}_1,E\tilde{x}_N\rangle - \min P]
\end{aligned}
\end{equation}
where the second inequality results from the convexity of $f(\cdot)$ and the fact $\hat{\eta}_1=\hat{\eta}_2=\cdots=\hat{\eta}_K$. The assertion \eqref{conv:sublinear} follows immediately from the above inequality. \endproof

Next, we shall prove the $o(1/\nu)$ convergence of the ADA.  Motivated by \cite{deng2017parallel,he2015non}, we use the quantity $||u^\nu-u^{\nu+1}||_G^2$ as a measure of the convergence rate. In fact, if $||u^\nu-u^{\nu+1}||_G^2=0$, then $u^{\nu+1}$ is an optimal solution, \ie, $(x^{\nu+1},\eta_1^{\nu+1})\in X^*\times Y^*$. More explicitly, $||u^\nu-u^{\nu+1}||_G^2=0$ implies the following:
\begin{equation}
x^\nu=x^{\nu+1}\text{ and } w^\nu=w^{\nu+1}.
\end{equation}
By the update step for $w$, we can conclude $\eta_1^{\nu+1}=\cdots=\eta_K^{\nu+1}$. Combining this with $x^\nu=x^{\nu+1}$, we obtain
\begin{equation}
0\in\partial f(x^{\nu+1})+E^T\eta_1^{\nu+1}+N_{X}(x^{\nu+1}),
\end{equation}
or equivalently, $(x^{\nu+1},\eta_1^{\nu+1})\in X^*\times Y^*$. Conversely, if the quantity $||u^\nu-u^{\nu+1}||_G^2$ is relatively large, $u^{\nu+1}$ should not be close to the optimal solution set. Based on previous analysis, $||u^\nu-u^{\nu+1}||_G^2$ is a reasonable measure to quantify the distance between $u^{\nu+1}$ and the optimal solution set.

To show the convergence rate, we first prove the following lemma on the monotonicity property of the iterations:
\begin{lemma}\label{lem:mono_u}
	Let $u^\nu$ be defined as in Theorem \ref{thm:sublinear}. Then
	\be\label{mono_u}
	||u^\nu-u^{\nu+1}||_G^2\le||u^{\nu-1}-u^{\nu}||_G^2.\ee
\end{lemma}
\proof For notational simplicity, for each iteration $\nu$, we introduce 
\begin{equation}\label{diff:u}
\Delta u^{\nu+1}=\begin{pmatrix}\Delta w^{\nu+1}\\ \Delta x^{\nu+1}\\ \Delta \eta^{\nu+1}\\ \Delta \zeta^{\nu+1}\end{pmatrix}=\begin{pmatrix}w^\nu- w^{\nu+1}\\ x^\nu- x^{\nu+1}\\ \eta^\nu- \eta^{\nu+1}\\ \zeta^\nu- \zeta^{\nu+1} \end{pmatrix}. 
\end{equation}
By the optimality of $x_k^{\nu+1}$ in iteration $\nu+1$ and the update rule of $\eta_k^{\nu+1}$, we have
\begin{equation}\label{opt:x_k}
\frac{1}{c}(x_k^{\nu}-x_k^{\nu+1})-E_k^T\eta_k^{\nu+1}\in\partial{f_k(x_k^{\nu+1})}+N_{X_k}(x_k^{\nu+1})
,\qquad k=1,\ldots,K.
\end{equation}
Considering the $\nu$-th and $\nu+1$-th iteration, such optimality yields
\begin{equation}\label{opt:twoiters}
\underbrace{\frac{1}{c}\langle \Delta x_k^{\nu+1},\Delta x_k^\nu-\Delta x_k^{\nu+1}\rangle}_{(a)}-\langle E_k\Delta x_k^{\nu+1},\Delta \eta_k^{\nu+1}\rangle\ge 0,\qquad k=1,\ldots,K.
\end{equation}
For the second term in the above inequality,
\begin{equation}
\begin{aligned}
&-\sum_{k=1}^{K}\langle E_k\Delta x_k^{\nu+1},\Delta \eta_k^{\nu+1}\rangle=\sum_{k=1}^{K}\langle E_kx_k^{\nu+1}-E_kx_k^{\nu},\Delta \eta_k^{\nu+1}\rangle\\=&
\sum_{k=1}^{K}\langle \frac{\eta_k^{\nu+1}-y_k^\nu}{\rho/2}-\frac{\eta_k^{\nu}-y_k^{\nu-1}}{\rho/2}+w_k^\nu-w_k^{\nu-1},\Delta \eta_k^{\nu+1}\rangle\\=&
\sum_{k=1}^{K}\langle \frac{\eta_k^{\nu+1}-\frac{\eta_k^\nu+\zeta_k^\nu}{2}}{\rho/2}-\frac{\eta_k^{\nu}-\frac{\eta_k^{\nu-1}+\zeta_k^{\nu-1}}{2}}{\rho/2}+\frac{\eta_k^\nu-\zeta_k^\nu}{\rho},\Delta \eta_k^{\nu+1}\rangle\\=&\sum_{k=1}^{K}\frac{1}{\rho}\langle\Delta\eta_k^\nu-\Delta\eta_k^{\nu+1},\Delta \eta_k^{\nu+1}\rangle+\sum_{k=1}^{K}\frac{1}{\rho}\langle\Delta\zeta_k^\nu,\Delta\eta_k^{\nu+1}\rangle+\\&\sum_{k=1}^{K}\frac{1}{\rho}\langle\eta_k^{\nu+1}-\eta_k^\nu,\eta_k^\nu-\eta_k^{\nu+1}\rangle+\sum_{k=1}^{K}\frac{1}{\rho}\langle\eta_k^\nu-\zeta_k^\nu,\Delta \eta_k^{\nu+1}\rangle\\=&\underbrace{\sum_{k=1}^{K}\frac{1}{\rho}\langle\Delta\eta_k^\nu-\Delta\eta_k^{\nu+1},\Delta \eta_k^{\nu+1}\rangle}_{(b)}+\underbrace{\sum_{k=1}^{K}\frac{1}{\rho}\langle\Delta\zeta_k^\nu-\Delta\zeta_k^{\nu+1},\Delta\zeta_k^{\nu+1}\rangle}_{(c)}+\\&\underbrace{\sum_{k=1}^{K}\frac{1}{\rho}\langle\eta_k^{\nu+1}-\zeta_k^{\nu+1},\eta_k^\nu-\eta_k^{\nu+1}\rangle}_{(d)}.
\end{aligned}
\end{equation}
For term (d),
\begin{equation}
\begin{aligned}
&2\sum_{k=1}^{K}\frac{1}{\rho}\langle\eta_k^{\nu+1}-\zeta_k^{\nu+1},\eta_k^\nu-\eta_k^{\nu+1}\rangle\\=&\sum_{k=1}^{K}\frac{1}{\rho}||\eta_k^\nu-\zeta_k^{\nu+1}||_2^2-\sum_{k=1}^{K}\frac{1}{\rho}||\eta_k^\nu-\eta_k^{\nu+1}||_2^2-\sum_{k=1}^{K}\frac{1}{\rho}||\eta_k^{\nu+1}-\zeta_k^{\nu+1}||_2^2\\=&\sum_{k=1}^{K}\frac{1}{\rho}||\eta_k^\nu-\zeta_k^\nu+\zeta_k^\nu-\zeta_k^{\nu+1}||_2^2-\sum_{k=1}^{K}\frac{1}{\rho}||\eta_k^\nu-\eta_k^{\nu+1}||_2^2-\sum_{k=1}^{K}\frac{1}{\rho}||\eta_k^{\nu+1}-\zeta_k^{\nu+1}||_2^2\\=&\sum_{k=1}^{K}\frac{1}{\rho}||\eta_k^\nu-\zeta_k^\nu||_2^2+\underbrace{\sum_{k=1}^{K}\frac{2}{\rho}\langle\eta_k^\nu-\zeta_k^\nu,\zeta_k^\nu-\zeta_k^{\nu+1}\rangle}_{=0}+\frac{1}{\rho}||\zeta^\nu-\zeta^{\nu+1}||_2^2-\\&\sum_{k=1}^{K}\frac{1}{\rho}||\eta_k^\nu-\eta_k^{\nu+1}||_2^2-\sum_{k=1}^{K}\frac{1}{\rho}||\eta_k^{\nu+1}-\zeta_k^{\nu+1}||_2^2.
\end{aligned}
\end{equation}
Applying the equality \eqref{3identity} to (a), (b) and (c) and combining them with the above transformation for term (d), the inequality \eqref{opt:twoiters} yields
\begin{equation}\label{G-norm:monotoncity}
\begin{aligned}
&||\Delta u^\nu||_G^2-||\Delta u^{\nu+1}||_G^2\ge\sum_{k=1}^{K}\frac{1}{c}||\Delta x_k^\nu-\Delta x_k^{\nu+1}||_2^2+\sum_{k=1}^{K}\frac{1}{\rho}||\Delta\eta_k^\nu-\Delta\eta_k^{\nu+1}||_2^2+\\&
\sum_{k=1}^{K}\frac{1}{\rho}||\Delta\zeta^\nu-\Delta\zeta^{\nu+1}||_2^2+\underbrace{\sum_{k=1}^{K}\frac{1}{\rho}(||\eta_k^\nu-\eta_k^{\nu+1}||_2^2-||\zeta_k^\nu-\zeta_k^{\nu+1}||_2^2)}_{\ge 0}\ge 0.
\end{aligned}
\end{equation}
The nonnegativity of the last term is a direct result of the definition  $\zeta_1^\nu=\cdots=\zeta_K^\nu=\frac{1}{K}\sum_{j=1}^{K}\eta_{j}^{\nu}$ and Cauchy–Schwarz inequality.  Hence the inequality \eqref{mono_u} holds.  \endproof

The following elementary lemma helps to improve the convergence rate from $O(1/\nu)$ to $o(1/\nu)$.
\begin{lemma}\label{lem:o(1/k)}
	Suppose a sequence $\{a_\nu\}_{\nu=0}^\infty\subseteq \Rset$ satisfies the following: (a) $a_\nu\ge 0$; (b)
	$\sum_{\nu=0}^{\infty} a_\nu<\infty$;
	and (c) $a_\nu$ is monotonically non-increasing.
	Then, we have $a_\nu=o(1/\nu)$.
\end{lemma}
\proof See Lemma 1.1 in \cite{deng2017parallel}. \endproof

Combining the results from previous two lemmas, we present the $o(1/\nu)$ convergence of the ADA.
\begin{theorem}\label{thm:o(1/k)}
	Let $\{u^\nu=(w^\nu,x^\nu,\eta^\nu,\zeta^\nu)\}$ in $W\times X\times S$ be the infinite sequence generated by the ADA, then 
		\vspace{-2ex}
	\begin{equation}\label{o(1/k)}
	||u^\nu-u^{\nu+1}||_G^2=o(1/\nu)
		\vspace{-2ex}
	\end{equation}
	holds and thus
			\vspace{-2ex}
	\begin{equation}\label{o(1/k):x}
	||x^\nu-x^{\nu+1}||_2^2=o(1/\nu)
			\vspace{-2ex}
	\end{equation}
	and 
		\vspace{-2ex}
	\begin{equation}\label{o(1/k):constraint}
	||\sum_{k=1}^{K}E_kx_k^{\nu+1}||_2^2=o(1/\nu).
		\vspace{-2ex}
	\end{equation}
\end{theorem}
\proof
In the proof of Theorem \ref{thm:global_conv_ADA}, we have shown that 
	\vspace{-3ex}
\[
\sum_{\nu=0}^{\infty}||u^\nu-u^{\nu+1}||_G^2<\infty. 
	\vspace{-3ex}
\]On the other hand, Lemma \ref{lem:mono_u} proved the non-increasing property of $||u^\nu-u^{\nu+1}||_G^2$. Hence, \eqref{o(1/k)} follows directly from Lemma \ref{lem:o(1/k)} and then \eqref{o(1/k):x} holds. For the estimate for the constraint in \eqref{o(1/k):constraint}, we have
\begin{equation}
\begin{aligned}
&||\sum_{k=1}^{K}E_kx_k^{\nu+1}||_2^2=||\sum_{k=1}^{K}(E_kx_k^{\nu+1}-w_k^\nu)||_2^2=\frac{4}{\rho^2}||\sum_{k=1}^{K}(\eta_k^{\nu+1}-y_k^\nu)||_2^2\\\le&\frac{4K}{\rho^2}\sum_{k=1}^{K}||\eta_k^{\nu+1}-\frac{1}{2}(\eta_k^\nu+\zeta_k^\nu)||_2^2\\=&\frac{K}{\rho^2}\sum_{k=1}^{K}||\eta_k^{\nu+1}-\eta_k^\nu+\eta_k^{\nu+1}-\zeta_k^{\nu+1}+\zeta_k^{\nu+1}-\zeta_k^\nu||_2^2\\\le&\frac{3K}{\rho^2}||\eta^{\nu+1}-\eta^\nu||_2^2+3K||w^{\nu+1}-w^\nu||_2^2+\frac{3K}{\rho^2}||\zeta^{\nu+1}-\zeta^\nu||_2^2=o(1/\nu),
\end{aligned}
\end{equation}  where the first two equalities result from $w^\nu\in W$ and the updating rule for $\eta^{\nu+1}$. This finishes the proof for \eqref{o(1/k):constraint}.\endproof

From Theorem \ref{thm:o(1/k)}, a reasonable stopping criterion for the ADA can be either
\begin{equation}\label{stop_criterion_A}
\frac{||x^\nu-x^{\nu+1}||}{\max\{1,||x^\nu||\}}\le \epsilon
\end{equation}
or
\begin{equation}\label{stop_criterion_B}
\frac{||Ex^{\nu+1}-q||}{\max\{1,||q||\}}\le \epsilon
\end{equation}
for some given tolerance $\epsilon$.

\subsection{Relation to the ADMM}\label{Relation: ADMM}
The ADA is closely related to the ADMM. Here, we compare the ADA with two variants of ADMM, namely, the Variable Splitting ADMM and the Proximal Jacobian ADMM. For simplicity of notation, we assume $q=0$.

Applying the classical two-block ADMM to the transformation in \eqref{P:w}, \cite{wang2013solving} proposed the following Variable Splitting ADMM (VSADMM), see Algorithm \ref{VSADMM}. 
\begin{algorithm}
	\caption{Variable Splitting ADMM}\label{VSADMM}
	\begin{algorithmic}[1]
		\State Given $w^0\in W, x^0\in X, y^0\in (\Rset^m)^K, \beta>0$
		\For{$\nu=0,1, \dots $}
		\State $x_k^{\nu+1} = \argmin_{x_k\in X_k}
		f_k(x_k)+\frac{\beta}{2}||E_kx_k-w_k^\nu+\frac{y_k^\nu}{\beta}||_2^2,\quad k=1,\ldots,K,$
		\State $w^{\nu+1}=\argmin_{w\in W}\frac{\beta}{2}\sum_{k=1}^{K}||E_kx_k^{\nu+1}-w_k+\frac{y_k^\nu}{\beta}||_2^2,$
		\State $y_k^{\nu+1}=y_k^{\nu}+\beta[E_kx_k^{\nu+1}-w_k^{\nu+1}],\quad k=1,\ldots,K.$
		\EndFor
	\end{algorithmic}
\end{algorithm}
The convergence result for VSADMM was established  on the basis of the classical two-block ADMM. Compared to the ADA, we notice that no proximal terms exist during the $x$-update in the VSADMM. Therefore, the full column rank assumption of $E_k$ is necessary for the VSADMM to guarantee the solution uniqueness in each iteration. The $w$-update step in the VSADMM also differs from that in the ADA as it does not use the information on the previous iteration explicitly. 

The Proximal Jacobian ADMM (Prox-JADMM) provided in \cite{deng2017parallel} solves problem \eqref{P} directly by adding a proximal term in the Jacobian-type ADMM, see Algorithm \ref{PJADMM}. 
\begin{algorithm}
	\caption{Proximal Jacobian ADMM}\label{PJADMM}
	\begin{algorithmic}[1]
		\State Given $ x^0\in X, \lambda^0\in\Rset^m, \beta>0$
		\For{$\nu=0,1, \dots $}
		\For  {$k=1,\ldots,K$}
		\State $x_k^{\nu+1} = \argmin_{x_k\in X_k}
		f_k(x_k)+\frac{\beta}{2}||E_kx_k+\sum_{j\neq k}E_jx_j^\nu-\frac{\lambda^\nu}{\beta}||_2^2+
		\frac{1}{2}||x_k-x_k^\nu||_{P_k}^2,$
		\EndFor
		\State $\lambda^{\nu+1}=\lambda^{\nu}-\gamma\beta\sum_{k=1}^{K}E_kx_k^{\nu+1},$
		\EndFor
	\end{algorithmic}
\end{algorithm}
It is worth noting that the ADA shares the same $o(1/\nu)$ convergence rate  as the Prox-JADMM.  
However, the Prox-JADMM requires
the constraints $E_k$, the proximal terms $P_k$  and the damping parameter $\gamma$ to satisfy certain relationships to guarantee 
the convergence. Because the convergence results for the ADA are established using a very different approach,
we impose no restriction on the proximal terms.

\section{The Inexact Augmented Decomposition Algorithm}\label{sec:iADA}
Here, we first review the general convergence theory of the (inexact-)proximal point algorithm (PPA) developed in \cite{rockafellar1976augmented,rockafellar1976monotone}. Let  $\mathcal{T}:\mathcal{X}\rightrightarrows \mathcal{X}$ be a maximally monotone operator. In order to solve the inclusion problem:
\begin{equation}\label{0inT(x)}
0\in \mathcal{T}(z),
\end{equation}
PPA takes the form of 
\begin{equation}\label{PPAstep}
z^{k+1}\approx (I+c_k\mathcal{T})^{-1}(z^k), \quad \forall k\ge 0,
\end{equation}
in the $(k+1)$-th iteration with a given sequence $c_k\uparrow c_\infty \le \infty$. The convergence result of PPA can be guaranteed as long as the approximation computation satisfies certain criteria; see \cite{rockafellar1976augmented,rockafellar1976monotone}. In addition, the local linear convergence result could be established when $\mathcal{T}^{-1}$ is Lipschitz continuous at the origin. In accordance with the PPA, the inexact version of the ADA comes out naturally as follows in Algorithm \ref{iexactADA}. The iADA allows the subproblems to be solved inexactly which is very important in many applications as it might be very expensive to solve these subproblems exactly.
\begin{algorithm}
	\caption{Inexact augmented decomposition algorithm}\label{iexactADA}
	\begin{algorithmic}[1]
		\State Given $w^0\in W, x^0\in X, y^0\in (\Rset^m)^K$
		\For{$\nu=0,1, \dots $}
		\State $x_k^{\nu+1} \approx \argmin_{x_k\in X_k}\phi_{k,\rho,c}^\nu(x_k),  k=1,\ldots,K$
		\State $\eta_{k}^{\nu+1}=\begin{cases}
		y_{k}^\nu+\frac{\rho}{2}[E_kx_k^{\nu+1}-w_{k}^\nu],\text{ if }k=1,\ldots,K-1\\
		y_{K}^\nu+\frac{\rho}{2}[E_Kx_K^{\nu+1}-q-w_{K}^\nu],\text{ if }k=K
		\end{cases}$
		\For  {$k=1,\ldots,K$}
		\State $\zeta_k^{\nu+1}= \frac{1}{K}\sum_{j=1}^{K}\eta_{j}^{\nu+1}$\\
		\State$w_{k}^{\nu+1}=w_{k}^\nu+\frac{1}{\rho}[\eta_{k}^{\nu+1}-\zeta_k^{\nu+1}]$\\
		\State$y_{k}^{\nu+1}=\frac{1}{2}[\eta_{k}^{\nu+1}+\zeta_k^{\nu+1}]$
		\EndFor
		\EndFor
	\end{algorithmic}
\end{algorithm}

Two natural concerns arise for the iADA: (1) the global convergence and (2) the local convergence rate. For that, we make the following assumptions on $f$ for the rest of the paper: 
\begin{assumption}\label{assu1}
	(a)  $f=f_1(x_1)+\cdots+f_K(x_K)$, with each $f_k$ given by
	\begin{equation}
	f_k(x_k)=g_k(A_kx_k)+h_k(x_k)
	\end{equation}
	where $g_k$ and $h_k$ are both closed proper convex functions and $A_k$'s are some given matrices.\\
	(b) Every $g_k$ is strongly convex and continuously differentiable on \intn\textnormal{(\dom$g_k$)} with a Lipschitz continuous gradient
	\begin{equation}
	||A_k^T\nabla g_k(A_kx_k)-A_k^T\nabla g_k(A_kx_k')||\le L_g^k||A_k(x_k-x_k')||,\qquad \forall x_k,x_k'\in X_k
	\end{equation} 
	where $L_g^k\ge 0,k=1,\ldots,K$.\\
	(c) The epigraph of each $h_k$ is a polyhedral convex set.\\
	(d) The feasible sets $X_k,k=1,\ldots,K$ are polyhedral convex sets.\\
	(e) The feasible sets $X_k,k=1,\ldots,K$ are compact sets.
\end{assumption}
Here are several comments on the above assumptions.
\begin{itemize}
	\item Either  $g_k$ or $h_k$ can be absent in $f_k$. Although $g_k$ is assumed to be strongly convex, we do not impose any condition on $A_k$. Therefore, $f_k$ is not necessarily strongly convex in general and the optimal solution is not necessarily unique. 
	\item We do not assume any condition for the rank of $E_k,k=1,\ldots,K$ which is required to have full column rank in \cite{hong2017linear}. For the ADMM, this assumption is necessary to ensure that in each iteration, the subproblem for the $k$-th block is strongly convex. But for the iADA, this assumption is no longer required as there exists a proximal term in each subproblem which makes its optimality attainable and unique.
	\item The compactness assumption of $X_k,k=1,\ldots,K$ will facilitate the proof in Section \ref{sec:stability} and is not necessary for the convergence result in Section \ref{sec:Conv_iADA} due to the boundedness of the sequence generated by the iADA.
\end{itemize}   Based on these assumptions, we can simply write $f$ as
\begin{equation}\label{f_A}
f(x)=g(Ax)+h(x)=\sum_{k=1}^{K}g_k(A_kx_k)+\sum_{k=1}^{K}h_k(x_k)
\end{equation}
where $g(Ax)=\sum_{k=1}^{K}g_k(A_kx_k)$ and $h(x)=\sum_{k=1}^{K}h_k(x_k)$ represent the smooth and nonsmooth parts respectively. In addition, $g(\cdot)$ is strongly convex and $h(\cdot)$ is convex with a polyhedral epigraph. The strong convexity of $g(\cdot)$ implies the following proposition, whose proof is omitted.
\begin{proposition}\label{prop:x_constant}
	For any $x$ in the solution set $X^*$, $A_kx_k,k=1,\ldots,K$ are constant and hence $Ax$ is constant. 
\end{proposition}

In the next section, we will discuss the stability result of the Lagrangian function under some perturbations which is essential to the local linear convergence result. 

\section{On the stability results of $\mathcal{T}_{\bar{L}}$}\label{sec:stability}
In this section, we establish the stability result of the maximal monotone operator $\mathcal{T}_{\bar{L}}$ defined in \eqref{T_bar_L} corresponding to the perturbations of both primal and dual solutions under Assumption \ref{assu1}. This property serves the key ingredient for the local convergence rate analysis of the iADA. 

Recall the definition of $\bar{L}(w,x,\eta,\zeta)$ in \eqref{L:final_v}. For each $(w,x,\eta,\zeta)\in W\times X\times S$, $\mathcal{T}_{\bar{L}}(w,x,\eta,\zeta)$ is defined as
\vspace{-1ex}
\begin{equation}\label{T_bar_L_equiv}
\vspace{-1ex}
T_{\bar{L}}(w,x,\eta,\zeta)=\{(v_1,v_2,v_3,v_4)|(v_1,v_2,-v_3,-v_4)\in \partial \bar{L}(w,x,\eta,\zeta)\},
\end{equation}
or equivalently,  $\mathcal{T}_{\bar{L}}(w,x,\eta,\zeta)$ is  the set of $v=(v_1,v_2,v_3,v_4)\in(\Rset^m)^K\times\Rset^n\times(\Rset^m)^K\times(\Rset^m)^K$ such that 
\vspace{-1ex}
\begin{equation}\label{T_bar_L}
\begin{aligned}
&\bar{L}(w',x',\eta,\zeta)-\langle w',v_1\rangle-\langle x',v_2\rangle+\langle \eta,v_3\rangle+\langle \zeta,v_4\rangle\\
\ge &\bar{L}(w,x,\eta,\zeta)-\langle w,v_1\rangle-\langle x,v_2\rangle+\langle \eta,v_3\rangle+\langle \zeta,v_4\rangle\\
\ge & \bar{L}(w,x,\eta',\zeta')-\langle w,v_1\rangle-\langle x,v_2\rangle+\langle \eta',v_3\rangle+\langle \zeta',v_4\rangle \\  &\text{ for all } (w',x')\in W\times X,(\eta',\zeta')\in S.
\end{aligned}
\end{equation}
Any solution to $(0,0,0,0)\in \mathcal{T}_{\bar{L}}(w,x,\eta,\zeta)$ is a saddle point of $\bar{L}$. Denote $v_1=(v_{1,1},\ldots,v_{1,K})\in(\Rset^m)^K, v_2=(v_{2,1},\ldots,v_{2,K})\in\Rset^{n_1}\times\cdots\times\Rset^{n_K}$,  $v_3=(v_{3,1},\ldots,v_{3,K})\in(\Rset^m)^K$ and $P_{W^\perp}(v_4)=(v_4^\perp,\ldots,v_4^\perp)\in(\Rset^m)^K$. We consider the following perturbed form of problem \eqref{P:w}:
\begin{equation}\label{P:w_p}
\begin{aligned}
& \underset{x,w}{\text{min}}
& & f_{1}(x_1)+\cdots+f_{K}(x_K)-\langle w,v_1\rangle-\langle x,v_2\rangle\\
& \text{ s.t.}
& & E_kx_k-w_k +v_{3,k}+v_4^\perp= 0 ,\quad k = 1, \ldots, K-1,\\
&
& & E_Kx_K-q-w_K+v_{3,K}+v_4^\perp= 0,\\
&  
& &w_1+\cdots+w_K = 0,\\
&
& & x_k\in X_k,\quad k=1,2,\ldots,K
\end{aligned}
\end{equation}
Its corresponding KKT conditions are given by 
\begin{equation}\label{KKT:w_p}
\begin{aligned}
-E_k^T\eta_k+v_{2,k}\in\partial f_k(x_k)+N_{X_k}(x_k), &\quad k=1,\ldots,K\\
-\eta_k+\mu=v_{1,k},& \quad k=1,\ldots,K\\
E_kx_k-w_k +v_{3,k}+v_4^\perp= 0,& \quad k=1,\ldots,K-1\\
E_Kx_K-q-w_K+v_{3,K}+v_4^\perp= 0,&\\
w_1+\cdots+w_K = 0.&
\end{aligned}
\end{equation}
One can easily check that 
\begin{equation}\label{T_bar_L_KKT}
\begin{aligned}
\mathcal{T}_{\bar{L}}^{-1}(v_1,v_2,v_3,v_4) = \text{ set of all } (w,x,\eta,P_{W^\perp}(\eta))\in{W\times X\times S} \\ \text{ such that there exists }  \mu\in\Rset^m  \text{ satisfying that }(w,x,\eta,\mu) \\ \text{ is a solution of the KKT conditions \eqref{KKT:w_p}}.
\end{aligned}
\end{equation}

Based on the above observation, we first study the stability results of the KKT system \eqref{KKT:w_p} under perturbations considered above. Under Assumption \ref{assu1}, every $f_k(x_k)$ is the sum of a smooth function $g_k(A_kx_k)$ and a nonsmooth function $h_k(x_k)$ with a polyhedral epigraph. By introducing a variable $s=(s_1,\ldots,s_K)\in\Rset^K$, for each $k$, we can rewrite the polyhedral set $\{(x_k,s_k):x_k\in X_k, h_k(x_k)\le s_k\}$ compactly as $C_x^kx_k+C_s^ks_k\ge c_k$ for some matrices $C_x^k\in\Rset^{j_k\times n_k}$, $C_s^k\in\Rset^{j_k\times 1}$ and $c_k\in\Rset^{j_k\times 1}$, where $j_k$s are some positive integers with $\sum_{k=1}^K j_k=j$. Then, we can transform \eqref{P:w} equivalently into 
\begin{equation}\label{P:w_eqv}
\begin{aligned}
& \underset{x,w,s}{\text{min}}
& & \sum_{k=1}^{K}g_k(A_kx_k)+s_k\\
& \text{ s.t.}
& & E_kx_k-w_k = 0 ,\quad k = 1, \ldots, K-1,\\
&
& & E_Kx_K-q-w_K = 0,\\
&  
& &w_1+\cdots+w_K = 0,\\
&
& & C_x^kx_k+C_s^ks_k-c_k\ge 0,\quad k=1,2,\ldots,K.
\end{aligned}
\end{equation}
For the perturbed problem \eqref{P:w_p}, similarly, we have the following equivalent transformation:
\begin{equation}\label{P:w_p_eqv}
\begin{aligned}
& \underset{x,w,s}{\text{min}}
& & \sum_{k=1}^{K}g_k(A_kx_k)+s_k-\langle w,v_1\rangle-\langle x,v_2\rangle\\
& \text{ s.t.}
& & E_kx_k-w_k +v_{3,k}+v_4^\perp= 0 ,\quad k = 1, \ldots, K-1,\\
&
& & E_Kx_K-q-w_K +v_{3,K}+v_4^\perp= 0,\\
&  
& &w_1+\cdots+w_K = 0,\\
&
& & C_x^kx_k+C_s^ks_k-c_k\ge 0,\quad k=1,2,\ldots,K.
\end{aligned}
\end{equation}
The canonical Lagrangian function for \eqref{P:w_p_eqv} is
given by 
\be\label{L:w_p_eqv}
\begin{aligned}
	&L^v(w,x,s,\eta,\lambda,\mu)=\sum_{k=1}^{K}g_k(A_kx_k)+s_k-\langle w,v_1\rangle-\langle x,v_2\rangle \\+&\sum_{k=1}^{K-1}\langle E_kx_k-w_k+v_{3,k}+v_4^\perp,\eta_k \rangle+ \langle E_Kx_K-q-w_K +v_{3,K}+v_4^\perp, \eta_K\rangle\\-&\sum_{k=1}^K\langle C_x^kx_k+C_s^ks_k-c_k,\lambda_k\rangle+\langle w_1+\cdots+w_K,\mu\rangle.
\end{aligned}
\ee
We use $Sol(P(v_1,v_2,v_3,v_4))$ to denote the set of saddle points for the Lagrangian function $L^v(w,x,s,\eta,\lambda,\mu)$  defined above corresponding to the perturbed problem \eqref{P:w_p_eqv}. Let $(v_1,v_2,v_3,v_4)=(0,0,0,0)$, then $Sol(P(0,0,0,0))$ represents the set of saddle points for the Lagrangian function of problem \eqref{P:w_eqv}. In order to show the stability results for the KKT system \eqref{KKT:w_p}, we define a set-valued mapping $\mathcal{M}$ that assigns the vector $(d,e,f)\in\Rset^n\times(\Rset^m)^K\times(\Rset^m)^K$ to the set of $(w,x,s,\eta,\lambda,\mu)\in(\Rset^m)^K\times\Rset^n\times\Rset^K\times(\Rset^m)^K\times\Rset^j\times\Rset^m$ that satisfy the following equations
\begin{equation}\label{M}
\begin{aligned}
-E_k^T\eta_k+(C_x^k)^T\lambda_k = d_k,& \quad k=1,\ldots,K\\
-\eta_k+\mu =e_k,& \quad k=1,\ldots,K\\
w_k-E_kx_k=f_k,& \quad k=1,\ldots,K-1\\
w_K-E_Kx_K+q=f_K,& \\
w_1+\cdots+w_K = 0,&\\
0\le \lambda_k\perp C_x^kx_k+C_s^ks_k-c_k\ge 0,& \quad k=1,\ldots,K\\
(C_s^k)^T\lambda_k =1,& \quad k=1,\ldots,K.
\end{aligned}
\end{equation}
One can easily verify that 
\begin{equation}
\begin{aligned}
(w,x,s,\eta,\lambda,\mu)\in\mathcal{M}(A^T\nabla g(Ax)-v_2,v_1,v_3+P_{W^\perp}(v_4)) \\
\text{if and only if } (w,x,s,\eta,\lambda,\mu)\in Sol(P(v_1,v_2,v_3,v_4)),
\end{aligned}
\end{equation}
\ie, a solution of the KKT system of \eqref{P:w_p_eqv} is also a saddle point of the Lagrangian function \eqref{L:w_p_eqv}. By taking $(v_1,v_2,v_3,v_4)=(0,0,0,0)$, we see that $(w^*,x^*,s^*,\eta^*,\lambda^*,\mu^*)\in\mathcal{M}(A^T\nabla g(Ax^*),0,0)$ if and only if $(w^*,x^*,s^*,\eta^*,\lambda^*,\mu^*)\in Sol(P(0,0,0,0))$. It is easily seen that $\mathcal{M}$ is a polyhedral multifunction; \ie, the graph of $\mathcal{M}$ is the union of a finitely many polyhedral convex sets. In \cite{robinson1981some}, Robinson established the following proposition that $\mathcal{M}$ enjoys the local upper Lipschitzian continuity property; see also \cite{hoffman2003approximate}.

\begin{proposition}\label{loc_up_lip}
	There exists a positive scalar $\theta$ that depends on $A,E,C_x,C_s$ only, such that for each $(\bar{d},\bar{e},\bar{f})$ there is a positive $\delta'$ satisfying 
	\begin{equation}
	\mathcal{M}(d,e,f)\subseteq\mathcal{M}(\bar{d},\bar{e},\bar{f})+\theta||(d,e,f)-(\bar{d},\bar{e},\bar{f})||\mathcal{B}\text{ whenever } ||(d,e,f)-(\bar{d},\bar{e},\bar{f})||\le\delta'
	\end{equation}
	where $\mathcal{B}$ is the unit Euclidean ball in $(\Rset^m)^K\times\Rset^n\times\Rset^k\times(\Rset^m)^K\times\Rset^j\times\Rset^m$.
\end{proposition}
Based on this proposition,  we claim that 
\begin{lemma}\label{KKT_eqv_lip}
	Suppose Assumptions \ref{assu:saddle-point} and \ref{assu1} hold. Then there exist positive scalars $\delta,\tau$ depending on $A,E,C_x,C_s$ only, such that for all $v=(v_1,v_2,v_3,v_4)\in(\Rset^m)^K\times\Rset^n\times(\Rset^m)^K\times(\Rset^m)^K$ and $||v||\le\delta$, any $(w(v),x(v),s(v),\eta(v),\lambda(v),\mu(v))\in Sol(P(v_1,v_2,v_3,v_4))$, we have
	 \begin{equation}
	dist((w(v),x(v),s(v),\eta(v),\lambda(v),\mu(v)),Sol(P(0,0,0,0)))\le\tau||v||.
	\end{equation}
\end{lemma}

\proof By the previous proposition, $\mathcal{M}$ is locally upper Lipschtizian with modulus $\theta$ at $(A^T\nabla g(Ax^*),0,0)$ for any $x^*\in X^*$. First we show that as $v\rightarrow 0$, $A^T\nabla g(Ax(v))\rightarrow A^T\nabla g(Ax^*)$. For that, take a sequence $v^i=(v_1^i,v_2^i,v_3^i,v_4^i)\in(\Rset^m)^K\times\Rset^n\times(\Rset^m)^K\times(\Rset^m)^K,i=1,2,\cdots,$ such that $||v^i||\rightarrow 0$. Based on Assumption \ref{assu1}(e), the sequence $x(v^i),i=1,2,\cdots$ lies in a compact set and so the other sequence $s(v^i)$ and $w(v^i)$ also belong to some compact sets, given the fact $s(v^i)=h(x(v^i))$ and the linear relationship among $x(v^i),v^i$ and $w(v^i)$. By passing to a subsequence if necessary, let $(w^\infty, x^\infty,s^\infty)$ be a cluster point of $\{(w(v^i),x(v^i),s(v^i)\}$. Due to the continuity of $\nabla g(\cdot)$, $(A^T\nabla g(Ax(v^i))-v_2^i,v_1^i,v_3^i+P_{W^\perp}(v_4^i))$ converges to $(A^T\nabla g(Ax^\infty),0,0)$ as $i\rightarrow \infty$. For all $i$, $\{(w(v^i),x(v^i),s(v^i),A^T\nabla g(Ax(v^i))-v_2^i,v_1^i,v_3^i+P_{W^\perp}(v_4^i))\}$ lies in the set $$\{(w,x,s,d,e,f)|(w,x,s,\eta,\lambda,\mu)\in \mathcal{M}(d,e,f) \text{ for some }(\eta,\lambda,\mu)\}$$ which is a closed polyhedral set. By passing to the limit, we can conclude $$(w^\infty,x^\infty,s^\infty,\eta^\infty,\lambda^\infty,\mu^\infty)\in \mathcal{M}(A^T\nabla g(Ax^\infty),0,0)$$
for some $(\eta^\infty,\lambda^\infty,\mu^\infty)\in (\Rset^m)^K\times\Rset^j\times\Rset^m$.
From Proposition \ref{prop:x_constant}, we know $Ax^\infty=Ax^*$ for any $x^*\in X^*$ which further implies that $A^T\nabla g(Ax(v))\rightarrow A^T\nabla g(Ax^*)$. Then there exists a positive scalar $\delta$ such that for all $v$ satisfying $||v||\le
\delta$, the following inequality $$||A^T\nabla g(Ax(v))-A^T\nabla g(Ax^*)||+||v||\le \delta' $$ holds. Based on Proposition \ref{loc_up_lip}, there exists $(w^*,x^*,s^*,\eta^*,\lambda^*,\mu^*)\in\mathcal{M}(A^T\nabla g(Ax^*),0,0)$, satisfying 
\be
\begin{aligned}
	||(w(v),x(v),s(v),\eta(v),\lambda(v),\mu(v))-(w^*,x^*,s^*,\eta^*,\lambda^*,\mu^*)||\\\le \theta(||A^T\nabla g(Ax(v))-A^T\nabla g(Ax^*)||+||v||).
\end{aligned}\ee
Since $(w(v),x(v),s(v),\eta(v),\lambda(v),\mu(v))\in\mathcal{M}(A^T\nabla g(Ax)-v_2,v_1,v_3+P_{W^\perp}(v_4))$, by the definition of $\mathcal{M}$ we have 
\begin{equation}\label{M_v}
\begin{aligned}
-E_k^T\eta_k(v)+(C_x^k)^T\lambda_k(v) = A_k^T\nabla g_k(A_kx_k(v))-v_2,& \quad k=1,\ldots,K\\
-\eta_k(v)+\mu(v) =v_{1,k},& \quad k=1,\ldots,K\\
w_k(v)-E_kx_k(v)=v_{3,k}+v_4^\perp,& \quad k=1,\ldots,K-1\\
w_K(v)-E_Kx_K(v)+q=v_{3,K}+v_4^\perp,& \\
w_1(v)+\cdots+w_K(v)= 0,&\\
0\le \lambda_k(v)\perp C_x^kx_k(v)+C_s^ks_k(v)-c_k\ge 0,& \quad k=1,\ldots,K\\
(C_s^k)^T\lambda_k(v) =1,& \quad k=1,\ldots,K.
\end{aligned}
\end{equation}
Similarly, since $(w^*,x^*,s^*,\eta^*,\lambda^*,\mu^*)\in\mathcal{M}(A^T\nabla g(Ax^*),0,0)$, it follows that 
\begin{equation}\label{M_*}
\begin{aligned}
-E_k^T\eta_k^*+(C_x^k)^T\lambda_k^* = A_k^T\nabla g_k(A_kx_k^*),& \quad k=1,\ldots,K\\
-\eta_k^*+\mu^* =0,& \quad k=1,\ldots,K\\
w_k^*-E_kx_k^*=0,& \quad k=1,\ldots,K-1\\
w_K^*-E_Kx_K^*+q=0,& \\
w_1^*+\cdots+w_K^* = 0,&\\
0\le \lambda_k^*\perp C_x^kx_k^*+C_s^ks_k^*-c_k\ge 0,& \quad k=1,\ldots,K\\
(C_s^k)^T\lambda_k^* =1,& \quad k=1,\ldots,K.
\end{aligned}
\end{equation}
Due to the strong convexity of  $g_k(\cdot)$ and the Lipschitzian continuity of its derivative $\nabla g_k(\cdot)$ in Assumption \ref{assu1}, there exist positive scalars $\sigma_g^k,L_g^k$ such that for all $x_1^k,x_2^k\in X_k$
\be\label{g_k_SC}
\langle A_k^T\nabla g_k(A_kx_1^k)-A_k^T\nabla g_k(A_kx_2^k),x_1^k-x_2^k\rangle\ge\sigma_g^k||A_kx_1^k-A_kx_2^k||^2,\ee and
\be\label{g_smooth}
||A_k^T\nabla g_k(A_kx_1^k)-A_k^T\nabla g_k(A_kx_2^k)||\le L_g^k||A_kx_1^k-A_kx_2^k||.
\ee 
Define $\sigma_g=\min_{k}\sigma_g^k$ and $L_g=\max_{k}L_g^k$.
Taking $x_1=x(v),x_2=x^*$, we obtain
\begin{equation*}\small
\begin{aligned}
&\sigma_g\sum_{k=1}^K||A_k(x(v)_k-x^*_k)||^2\\
\le &\sum_{k=1}^{K}\langle A_k^T\nabla g_k(A_kx(v)_k)-A_k^T\nabla g_k(A_kx^*_k),x(v)_k-x^*_k\rangle\\
=&\sum_{k=1}^{K}\langle -E_k^T(\eta(v)_k-\eta^*_k)+(C_x^k)^T(\lambda(v)_k-\lambda^*_k)+v_{2,k},x(v)_k-x^*_k\rangle\\
=&\sum_{k=1}^{K}\langle \lambda(v)_k-\lambda^*_k,C_x^kx(v)_k-C_x^kx^*_k\rangle+\sum_{k=1}^{K}\langle \eta(v)_k-\eta^*_k,-E_kx(v)_k+E_kx^*_k\rangle\\&+\sum_{k=1}^{K}\langle v_{2,k},x(v)_k-x^*_k\rangle
\end{aligned}
\end{equation*} 
where the first inequality comes from \eqref{g_k_SC} and the equalities come from \eqref{M_v} and \eqref{M_*}. Moreover, we have
\begin{equation*}
\begin{aligned}
&\sum_{k=1}^{K}\langle \lambda(v)_k-\lambda^*_k,C_x^kx(v)_k-C_x^kx^*_k\rangle\\
&=\sum_{k=1}^{K}\langle \lambda(v)_k-\lambda^*_k,C_x^kx(v)_k-C_x^kx^*_k\rangle+\langle \sum_{k=1}^{K}\lambda(v)_k-\lambda^*_k,C_s^ks(v)_k-C_s^ks^*_k\rangle\\
&=\sum_{k=1}^{K}\langle \lambda(v)_k-\lambda^*_k,(C_x^kx(v)_k+C_s^ks(v)_k-c_k)-(C_x^kx^*_k+C_s^ks^*_k-c_k)\rangle\\
&=-\sum_{k=1}^{K}[\langle \lambda^*_k,C_x^kx(v)_k+C_s^ks(v)_k-c_k\rangle+\langle \lambda(v)_k,C_x^kx^*_k+C_s^ks^*_k-c_k\rangle]\le 0
\end{aligned}
\end{equation*}
where the first equality follows from the fact that $(C_s^k)^T\lambda(v)_k=(C_s^k)^T\lambda^*_k=1,k=1,\ldots,K$ and the last equality and inequality both result from the complementary conditions in \eqref{M_v} and \eqref{M_*}. Consequently, we obtain that
\begin{equation*}\small
\begin{aligned}
&\sigma_g\sum_{k=1}^K||A_k(x(v)_k-x^*_k)||^2\\
&\le \sum_{k=1}^{K}\langle \eta(v)_k-\eta^*_k,-E_kx(v)_k+E_kx^*_k\rangle+\sum_{k=1}^{K}\langle v_{2,k},x(v)_k-x^*_k\rangle\\
&=\sum_{k=1}^{K}\langle \mu(v)+v_{1,k}-\mu^*,-w(v)_k+w^*_k+v_{3,k}+v_4^\perp\rangle+\sum_{k=1}^{K}\langle v_{2,k},x(v)_k-x^*_k\rangle\\ 
&= \sum_{k=1}^{K}\langle \mu(v)-\mu^*,v_{3,k}+v_4^\perp\rangle+\underbrace{\sum_{k=1}^{K}\langle \mu(v)-\mu^*,-w(v)_k+w^*_k\rangle}_{=0}\\&+\sum_{k=1}^{K}\langle v_{1,k},-w(v)_k+w^*_k+v_{3,k}+v_4^\perp\rangle+\sum_{k=1}^{K}\langle v_{2,k},x(v)_k-x^*_k\rangle\\
&\le||\mu(v)-\mu^*||(||v_3||+||v_4||)+||w(v)-w^*||||v_1||+||v_1||(||v_3||+||v_4||)+||(x(v)-x^*)||||v_2||\\
&\le||(w(v),x(v),\mu(v))-(w^*,x^*,\mu^*)||||v||+||v||^2.
\end{aligned}
\end{equation*}
Finally, based on Proposition \ref{loc_up_lip} and the above inequality, we have
\begin{equation*}
\begin{aligned}
&||(w(v),x(v),s(v),\eta(v),\lambda(v),\mu(v))-(w^*,x^*,s^*,\eta^*,\lambda^*,\mu^*)||^2\\&\le \theta^2(||A^T\nabla g(Ax(v))-A^T\nabla g(Ax^*)||+||v||)^2\\
&\le 2\theta^2(\sum_{k=1}^K||A_k^T\nabla g_k(A_kx(v)_k)-A_k^T\nabla g_k(A_kx_k^*)||^2+||v||^2)\\
&\le 2\theta^2(L_g^2\sum_{k=1}^K||A_k(x(v)_k-x^*_k)||^2+||v||^2)\\
&\le 2\theta^2 \max\{\frac{L_g^2}{\sigma_g},1\}(\sigma_g\sum_{k=1}^K||A_k(x(v)_k-x^*_k)||^2+||v||^2)\\
&\le 2\theta^2 \max\{\frac{L_g^2}{\sigma_g},1\}(||(w(v),x(v),\mu(v))-(w^*,x^*,\mu^*)||||v||+2||v||^2)\\
&\le 2\theta^2 \max\{\frac{L_g^2}{\sigma_g},1\}(||(w(v),x(v),s(v),\eta(v),\lambda(v),\mu(v))-(w^*,x^*,s^*,\eta^*,\lambda^*,\mu^*)||||v||+2||v||^2)
\end{aligned}.
\end{equation*}

We see the above inequality is quadratic in $ ||(w(v),x(v),s(v),\eta(v),\lambda(v),\mu(v))-(w^*,x^*,s^*,\eta^*,\lambda^*,\mu^*)||/||v||$, so we have
$$||(w(v),x(v),s(v),\eta(v),\lambda(v),\mu(v))-(w^*,x^*,s^*,\eta^*,\lambda^*,\mu^*)||/||v||\le \tau$$
for some scalar $\tau$ depending on $\theta,L_g,\sigma_g$ only. We conclude that $$dist((w(v),x(v),s(v),\eta(v),\lambda(v),\mu(v)),Sol(P(0,0,0,0)))\le\tau||v||.$$\endproof


In view of the operator $\mathcal{T}_{\bar{L}}$, combining  Lemma \ref{KKT_eqv_lip} with the observation in \eqref{T_bar_L_KKT}, we have the following corollary.

\begin{corollary}\label{cor:lip}
	Suppose Assumptions \ref{assu:saddle-point} and \ref{assu1} hold. Then there exist positive scalars $\delta,\tau$ depending on $A,E,C_x,C_s$ only, such that for all $v=(v_1,v_2,v_3,v_4)\in(\Rset^m)^K\times\Rset^n\times(\Rset^m)^K\times(\Rset^m)^K$ and $||v||\le\delta$, any $(w(v),x(v),\eta(v),\zeta(v))\in\mathcal{T}_{\bar{L}}^{-1}(v)$ satisfies
	\begin{equation}\label{T_bar_L_lip}
	dist(	(w(v),x(v),\eta(v),\zeta(v)),\mathcal{T}_{\bar{L}}^{-1}(0,0,0,0))\le 2\tau||v||.
	\end{equation}
\end{corollary}
\proof From Lemma \ref{KKT_eqv_lip} and observation in \eqref{T_bar_L_KKT} , we know that for any\\ $(w(v),x(v),\eta(v),\zeta(v))\in\mathcal{T}_{\bar{L}}^{-1}(v)$, there exists a $(w^*,x^*,\eta^*,\zeta^*)\in\mathcal{T}_{\bar{L}}^{-1}(0,0,0,0)$ satisfying that 
$$||(w(v),x(v),\eta(v))-(w^*,x^*,\eta^*)||\le \tau||v||.
$$
Since $\zeta(v)=P_{W^\perp}(\eta(v))$ and $\zeta^*=P_{W^\perp}(\eta^*)$, then
$$||(w(v),x(v),\eta(v),\zeta(v))-(w^*,x^*,\eta^*,\zeta^*)||\le 2\tau||v||
$$
holds which leads to \eqref{T_bar_L_lip}.\endproof

The compactness assumption of $X_k$ is indeed necessary for Corollary \ref{cor:lip}. However, if the generated sequence $\{x(v^i)\}$ lies in a compact set for a sequence $\{v^i\}_{i=1}^\infty$ converging to the origin, we claim the following result: under Assumptions \ref{assu:saddle-point} and \ref{assu1}(a)-(d), there exist positive scalars $\delta,\tau$ depending on $A,E,C_x,C_s$ only, when $||v^i||\le\delta$ the following
\be\label{T_bar_L_lip_PPA}
dist(	(w(v^i),x(v^i),\eta(v^i),\zeta(v^i)),\mathcal{T}_{\bar{L}}^{-1}(0,0,0,0))\le 2\tau||v^i||
\ee
holds. This observation relaxes the compactness assumption for $X_k,k=1,\ldots,K$ (Assumption \ref{assu1}(e)) when we show the local linear convergence in Theorem \ref{thm:loc_conv} for the iADA in Section \ref{sec:Conv_iADA}.

\section{Convergence analysis of the inexact ADA}\label{sec:Conv_iADA}

In this section, we study the convergence results of the inexact ADA for solving the problem \eqref{P}. For that, we first need to adopt the following stopping criterion developed in \cite{rockafellar1976augmented,rockafellar1976monotone} for approximately solving these subproblems
\begin{equation}\label{criteria_A}\tag{A}
\text{dist}(0,\partial\phi_{k,\rho,c}^\nu(x_k^{\nu+1}))\le \frac{\epsilon_\nu}{cK(\rho||E||+||E||+1)}, \qquad \sum_{\nu=0}^{\infty}\epsilon_\nu<\infty.
\end{equation}

\begin{theorem}\label{thm:global_conv}
	Suppose Assumption \ref{assu:saddle-point} holds and let $\{(w^\nu,x^\nu,\eta^\nu,\zeta^\nu)\}$ in $W\times X\times S$ be the infinite sequence generated by the ADA with the stopping criterion \eqref{criteria_A}. Then $(w^\nu,x^\nu,\eta^\nu,\zeta^\nu)$ converges to some saddle point $(\bar{w},\bar{x},\bar{\eta},\bar{\zeta})$ of \eqref{L:final_v} such that
	\begin{enumerate}[(a)]
		\item $(\bar{w},\bar{x})$ solves \eqref{P:w}, hence $\bar{x}$ solves \eqref{P},
		\item $\bar{\eta}_1=\cdots=\bar{\eta}_K\in \mathbb{R}^m$, and this common multiplier vector solves \eqref{D}. 	
	\end{enumerate}
\end{theorem}
\proof In each iteration $\nu$, we denote $(w_0^{\nu+1},x_0^{\nu+1},\eta_0^{\nu+1},\zeta_0^{\nu+1})=P_\nu(w^\nu,x^\nu,\eta^\nu,\zeta^\nu)$ as the exact saddle point of $\bar{L}^\nu(w,x,\eta,\zeta)$ and $(w^{\nu+1},x^{\nu+1},\eta^{\nu+1},\zeta^{\nu+1})$ as the inexact saddle point generated following the stopping criteria \eqref{criteria_A} respectively. By the update rule, the following estimates hold:
\begin{equation*}
||\eta_0^{\nu+1}-\eta^{\nu+1}||\le \frac{\rho||E||}{2}||x_0^{\nu+1}-x^{\nu+1}||,
\end{equation*} 
\begin{equation*}
||\zeta_0^{\nu+1}-\zeta^{\nu+1}||\le \frac{\rho||E||}{2}||x_0^{\nu+1}-x^{\nu+1}||,
\end{equation*}
and
\begin{equation*}
||w_0^{\nu+1}-w^{\nu+1}||\le ||E||||x_0^{\nu+1}-x^{\nu+1}||.
\end{equation*}
Thus, we can obtain
\begin{equation}\label{update_error}
||(w^{\nu+1},x^{\nu+1},\eta^{\nu+1},\zeta^{\nu+1})-P_\nu(w^\nu,x^\nu,\eta^\nu,\zeta^\nu)||\le (\rho||E||+||E||+1)||x^{\nu+1}-x_0^{\nu+1}||.
\end{equation}
Observing that the function $\phi_{k,\rho,c}^\nu$ defined in \eqref{phi_k} is strongly convex with modulus at least $\frac{1}{c}$ and $x_{0,k}^{\nu+1}$ minimize $\phi_{k,\rho,c}^\nu(x_k)$, we get 
\begin{equation}\label{SC:phi}
||x^{\nu+1}-x_0^{\nu+1}||\le c\sum_{k=1}^{K}\text{dist}(0,\partial\phi_{k,\rho,c}^{\nu+1}(x_k^{\nu+1})).
\end{equation}
Combining criterion \eqref{criteria_A}, \eqref{update_error} and \eqref{SC:phi}, we have 
\begin{equation}\label{update_error_final}
||(w^{\nu+1},x^{\nu+1},\eta^{\nu+1},\zeta^{\nu+1})-P_\nu(w^\nu,x^\nu,\eta^\nu,\zeta^\nu)||\le \epsilon_\nu, \text{ with } \sum_{\nu=1}^{\infty}\epsilon_\nu<\infty.
\end{equation}
From Assumption \ref{assu:saddle-point}, there exists a saddle point of the Lagrangian \eqref{L:P}. Therefore based on the relationship between \eqref{L:P} and \eqref{L:final_v} in Lemma \ref{L:relation}, there exists at least  one saddle point of the Lagrangian function $\bar{L}$. On the basis of \cite{rockafellar1976monotone}, the sequence of elements $(w^\nu,x^\nu,\eta^\nu,\zeta^\nu)$ generated in this manner  from any initial $(w^1,x^1)\in W\times X$ and $(\eta^1,\zeta^1)\in S$  converges to some saddle point $(\bar{w},\bar{x},\bar{\eta},\bar{\zeta})$ of the $\bar{L}$. Then $(\bar{w},\bar{x})$ solves \eqref{P:w} and  $(\bar{\eta},\bar{\zeta})$ solves \eqref{D:final_v}. By Lemma \ref{L:relation}, both (a) and (b) hold.\endproof

For the local convergence analysis, we need the following stopping criteria
\begin{equation}\label{criteria_B}\tag{B}
\text{dist}(0,\partial\phi_{k,\rho,c}^\nu(x_k^{\nu+1}))\le \frac{\epsilon'_\nu}{cK(\rho||E||+||E||+1)}\min\{1,||x_k^{\nu+1}-x_k^\nu||\}, \qquad \sum_{\nu=0}^{\infty}\epsilon'_\nu<\infty.
\end{equation}
The iADA does not impose any condition on the choice of $c$. We set $c=\rho$ for simplicity of the following analysis. The coefficient $\rho/2$ for the primal proximal term $||w-w^\nu||^2$ in \eqref{L:prox_v} can be changed to $1/2\rho$ after the rescaling $w'=\rho w$ and such rescaling only applies to the magnitude of $w$ and does not bring any other changes to the iADA. So this distinction from the standard proximal point method for minimax problems in \cite[Section 5]{rockafellar1976augmented} will not influence the following convergence results. 
\begin{theorem}\label{thm:loc_conv}
	Suppose Assumptions \ref{assu:saddle-point} and \ref{assu1} hold and let $\{(w^\nu,x^\nu,\eta^\nu,\zeta^\nu)\}$ in $W\times X\times S$ be the infinite sequence generated by the ADA with the stopping criterion \eqref{criteria_B}. Then, $(w^\nu,x^\nu,\eta^\nu,\zeta^\nu)$ converges to some saddle point $(\bar{w},\bar{x},\bar{\eta},\bar{\zeta})$ of \eqref{L:final_v} and  there exists $\{\theta_\nu\}$ such that
	$$dist((w^{\nu+1},x^{\nu+1},\eta^{\nu+1},\zeta^{\nu+1}),\mathcal{T}_{\bar{L}}^{-1}((0,0,0,0)))\le \theta_\nu dist((w^\nu,x^\nu,\eta^\nu,\zeta^\nu),\mathcal{T}_{\bar{L}}^{-1}((0,0,0,0)))$$
	for sufficient large $\nu$ and  $\lim_{\nu\rightarrow\infty}\theta_\nu=\frac{2\tau}{\sqrt{(4\tau^2+\rho^2)}}<1$ for some $\tau$.
\end{theorem}
\proof 
From Corollary \ref{T_bar_L_lip}, we have shown that there exist  $\tau,\delta>0$ such that for all $v=(v_1,v_2,v_3,v_4)\in(\Rset^m)^K\times\Rset^n\times(\Rset^m)^K\times(\Rset^m)^K$ and $||v||\le\delta$, any $(w(v),x(v),\eta(v),\zeta(v))\in\mathcal{T}_{\bar{L}}^{-1}(v)$ satisfies
\begin{equation}
dist((w(v),x(v),\eta(v),\zeta(v)),\mathcal{T}_{\bar{L}}^{-1}((0,0,0,0)))\le 2\tau||v||.
\end{equation}
So this theorem follows from \cite[Theorem 2.1]{luque1984asymptotic}.\endproof
\textbf{Remark 1.} In Theorem \ref{thm:global_conv}, we have shown that the sequence $\{u^\nu=(w^\nu,x^\nu,\eta^\nu,\zeta^\nu)\}$ converges to some saddle point $(\bar{w},\bar{x},\bar{\eta},\bar{\zeta})$ of \eqref{L:final_v} and hence $\{u^\nu\}$ lies in a compact set. Based on the observation in \eqref{T_bar_L_lip_PPA} and the proof of \cite[Theorem 2.1]{luque1984asymptotic}, the compactness of assumption of $X_k$ (Assumption \ref{assu1}(e)) is no longer needed for Theorem \ref{thm:loc_conv}.

\textbf{Remark 2.} When $c\neq\rho$, the local linear convergence still holds while the convergence rate ($\lim_{\nu\rightarrow\infty}\theta_\nu$) changes. 

Next, we provide some well-known examples on which the iADA enjoys the local linear convergence. 

%

\textbf{Convex regularization.} Many problems from empirical risk minimization and variable selection can be written as the following:
\begin{equation}\label{conv_reg}
\min_{x} f(x;(A,b))+r(x)
\end{equation}
where $A\in\Rset^{n\times d}$ and $b\in\Rset^n$, $f(\cdot)$ is the loss function which is often strongly convex with Lipschitz continuous gradient and $r(\cdot)$ is a convex regularization term which is possibly nonsmooth (\eg, the $\ell_1$-norm and TV-norm). By adding the constraint $x-z=0$, the above problem can be reformulate as 
\begin{equation}\label{conv_reg_v}
\begin{aligned}
& \min_{x,z} f(x;(A,b))+r(z)\\
& \text{ s.t.      }\quad  x-z=0.
\end{aligned}
\end{equation}

\textbf{Exchange problem.} Consider a network with $K$ agents exchanging $n$ commodities. Let $x_k\in\Rset^n$ be the amount of commodities in each agent $k$ and $f_k:\Rset^n\rightarrow \Rset$ be its corresponding cost function. The exchange problem is given by 
\be\label{exchange problem}
\begin{aligned}
	\min_{\{x_k\}_{k=1}^K}\sum_{k=1}^{K} f_k(x_k)\qquad \text{s.t. } \sum_{k=1}^{K}x_k=0
\end{aligned}
\ee
which minimizes the total cost  subject to the equilibrium constraint on all $K$ agents. In this special case, $E_k=\boldmath{I}$ and $q=0$. Optimization problems in this form arise in many areas such as resource allocation \cite{beck20141,xiao2006optimal}, multi-agent system \cite{you2011network} and image processing \cite{wright2012accelerated}. When the cost function $f_k$ in each agent satisfies Assumption \ref{assu1}(a)-(c), based on Theorem \ref{thm:loc_conv}, local linear convergence result is valid for the iADA under certain approximation criteria.

\section{Numerical Examples}\label{sec:numerical}
In this section, we demonstrate the linear convergence of both the exact ADA and the inexact ADA by some simple numerical examples. All the computational tasks for numerical experiments are implemented in Matlab 2017b running on a MacBook Pro. Retina, 2.6 GHz Intel Core i7 with 16Gb 2133 MHz LPDDR3 memory.
\subsection{The \textit{lasso} problem}
We perform some numerical experiments of Algorithm \ref{exactADA} for solving the following \textit{lasso} problem:
\vspace{-1.5ex}
\begin{equation}\label{lasso}
\min_{x\in\Rset^d} \frac{1}{2}||Ax-b||_2^2+\lambda_1||x||_1\\
\end{equation}
where  $A\in\Rset^{n\times d}$, $b\in\Rset^n$ and $\lambda_1$ is the regularization parameter. By introducing an auxiliary variable $z\in\Rset^d$, the above problem is equivalent to 
\vspace{-1.5ex}
\begin{equation}\label{lasso:v}
\vspace{-1.5ex}
\begin{aligned}
& \min_{x,z\in\Rset^d} \frac{1}{2}||Ax-b||_2^2+\lambda_1||z||_1\\
& \text{ s.t.      }\quad  x-z=0.
\end{aligned}
\end{equation}
Clearly, \eqref{lasso:v} is a two-block decomposition problem with $f_1(x_1)=\frac{1}{2}||Ax_1-b||_2^2$ and $f_2(x_2)=\lambda_1||x_2||_1$ by replacing $x \text{ and }z$ with $x_1 \text{ and }x_2$. Notice that $f_1$ and $f_2$ are not necessarily strongly convex. In this case,
\begin{equation}\label{lasso:phi_k}
\begin{aligned}
&\phi_{1,\rho,c}^\nu(x_1)=
\frac{1}{2}||Ax_1-b||_2^2+\frac{\rho}{4}||x_1-w_1^\nu+\frac{2}{\rho}y_1^\nu||_2^2+\frac{1}{2c}||x_1-x_1^\nu||_2^2 ,\\
&\phi_{2,\rho,c}^\nu(x_2)=\lambda_1||x_2||_1+\frac{\rho}{4}||x_2+w_2^\nu-\frac{2}{\rho}y_2^\nu||_2^2+\frac{1}{2c}||x_2-x_2^\nu||_2^2.
\end{aligned}
\end{equation}
For the first block, we can derive that 
\begin{equation}\label{lasso:phi_1}
x_1^{\nu+1}=[A^TA+(\frac{\rho}{2}+\frac{1}{c})\textbf{I}_d]^{-1}(A^Tb+\frac{\rho}{2}w_1^\nu+\frac{x_1^\nu}{c}-y_1^\nu).
\end{equation}
Though it may be time consuming to compute $[A^TA+(\frac{\rho}{2}+\frac{1}{c})\textbf{I}_d]^{-1}$ when $d$ is large, we only need to compute it at the initialization stage. The special structure of  $A^TA+(\frac{\rho}{2}+\frac{1}{c})\textbf{I}_d$ can be exploited and substantially improve performance, see \cite[Section 4.2]{boyd2011distributed}.
For the second block, the exact solution to the subproblem in each iteration is given by 
\begin{equation}\label{lasso:phi_2}
x_2^{\nu+1}:= S(\frac{y_2^\nu+x_2^\nu/c-\rho w_2^\nu/2}{\rho/2+1/c},\frac{\lambda_1}{\rho/2+1/c})
\end{equation}
where the \emph{soft thresholding operator} $S$ is defined in \cite{boyd2011distributed}.

 We generate the matrix $A$  and $0.05d$ nonzero entries of the sparse vector $x_0\in \Rset^d$ from the standard Gaussian distribution  $\mathcal{N}(0,1)$. We then let the response vector $b\in\Rset^{n}$ be given by $b=Ax_0+\epsilon$ where $\epsilon\sim \mathcal{N}(0,10^{-3}\textbf{I}_{n})$  and let the regularization parameter $\lambda_1$ be  $0.1||A^Tb||_\infty$.  We test the algorithm on two different sets of $(n,d)$: $(1000,4000)$, $(2000,20000)$.

\begin{figure}[h]
		\vspace{-3ex}
		\centering
		\includegraphics[width=.8\linewidth,height=0.25\textheight]{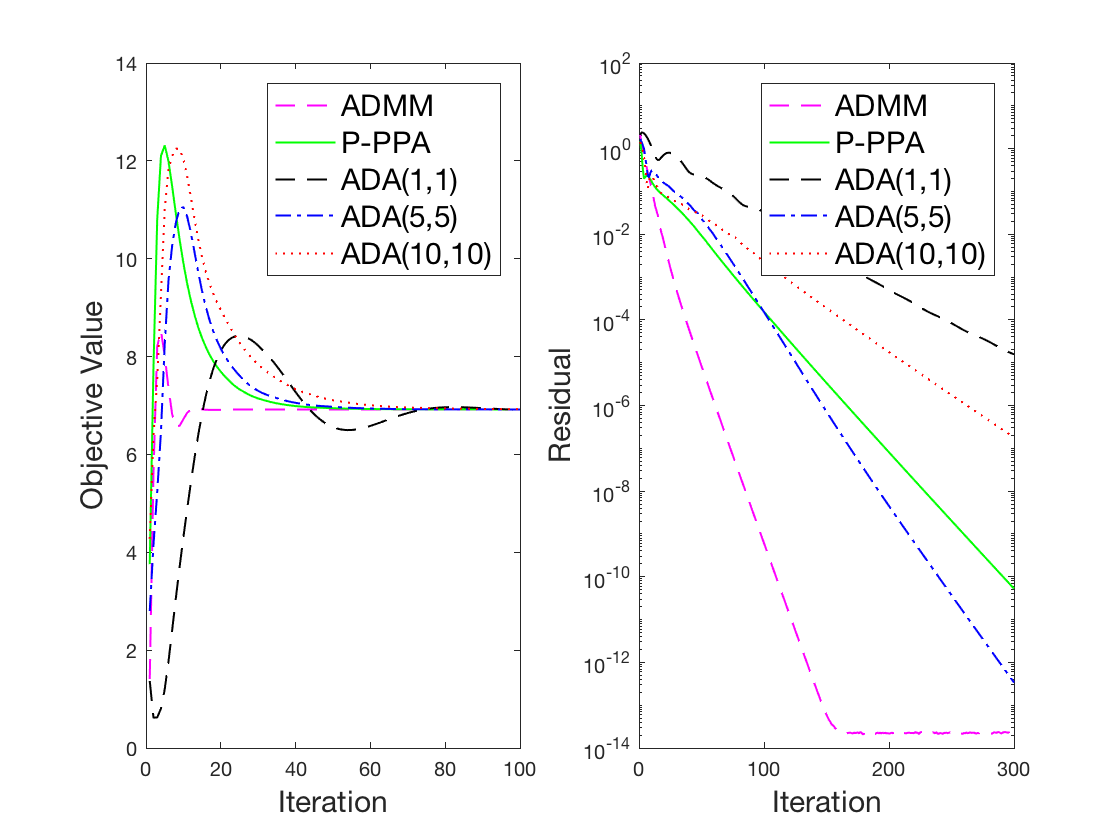}
		\caption{Convergence results of ADA, ADMM and P-PPA for the \textit{lasso}: $(n,d)=(1000,4000)$.}
		\label{fig1:sfig1}
\end{figure}

\begin{figure}[h]
		\centering
		\includegraphics[width=.8\linewidth,height=0.25\textheight]{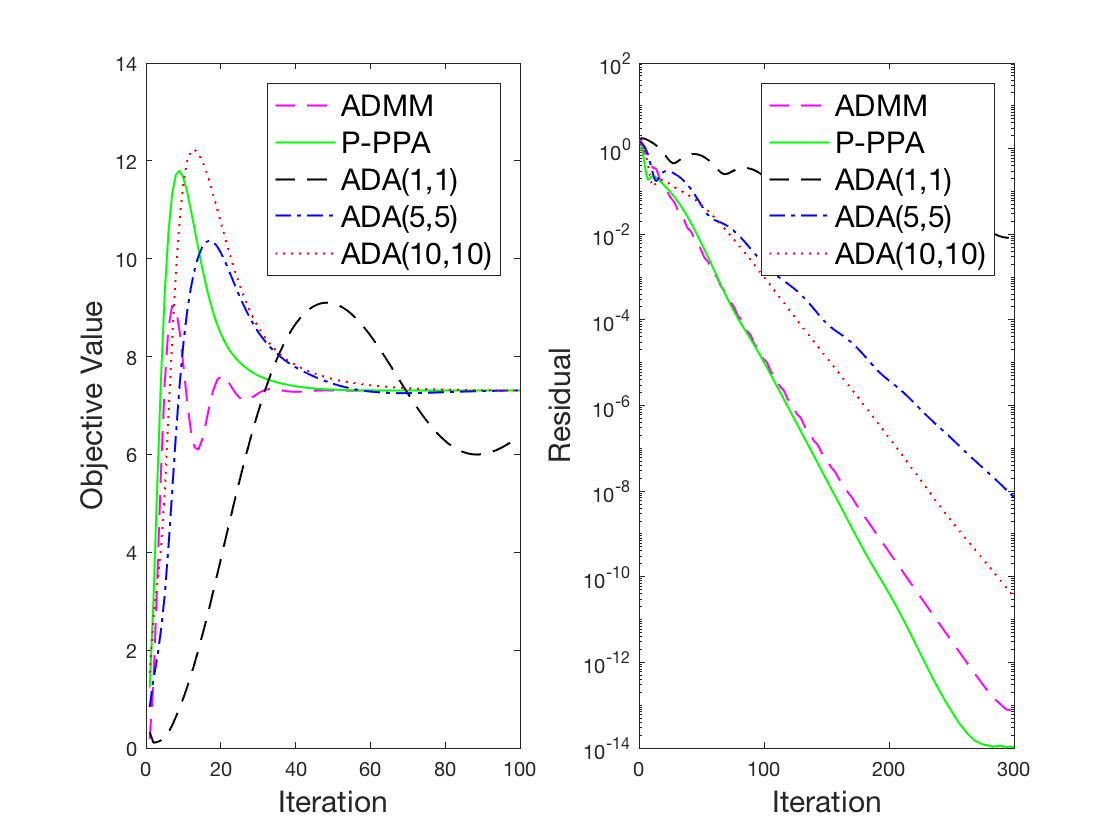}
		\caption{Convergence results of ADA, ADMM and P-PPA for the \textit{lasso}: $(n,d)=(2000,20000)$.}
		\label{fig1:sfig2}
\end{figure}

In our test, we compare the result of ADA with two other methods for the \textit{lasso} problem: ADMM\footnote[1]{Available at \url{http://web.stanford.edu/∼boyd/papers/admm/}}\cite{boyd2011distributed} and P-PPA\cite{bai2017parameterized}. For the implementation of  ADMM, we take a widely-used step-length 1.618 and a fixed penalty parameter 1. For P-PPA,  we used the parameters suggested in \cite{bai2017parameterized} for solving the \textit{lasso}. For the ADA, we choose the following three pairs of $(\rho, c)$:  $(1,1),(5,5),(10,10)$.  In each iteration, we solve both subproblems exactly and the computational time for all three algorithms is nearly the same.  For all algorithms, we use the same initial point $(x^0,y^0)=(\textbf{0},\textbf{0})$ and run 300 iterations. For all comparison algorithms, we report the objective value $f(x^\nu)=\frac{1}{2}||Ax_1^\nu-b||^2+\lambda_1||x_2^\nu||_1$,  and the residual norm $||x_1^\nu-x_2^\nu||$. The convergence results  are presented in Figures \ref{fig1:sfig1} and \ref{fig1:sfig2}.
	\vspace{-1ex}

	\vspace{-5ex}
From Figure \ref{fig1:sfig1}, we notice that ADMM performs best in the case $(n,d)=(1000,4000)$ while  ADA achieves comparable performance with P-PPA when $(\rho,c)=(5,5)$. This suggests that the convergence of ADA becomes slow if the proximal parameter is either too big or too small. When $(n,d)=(2000,20000)$, P-PPA shows the best convergence and ADA with $(\rho,c)=(10,10)$ converges a little bit slower. Both ADMM and P-PPA methods use the  Gauss-Seidel style update which tends to converge faster in terms of iterations, since it is able to incorporate information from the other coordinates more quickly. However, the Jacobi style update of  ADA is more amenable for parallelization. 

\subsection{The exchange problem} For the exchange problem in \eqref{exchange problem}, we consider the quadratic cost function \\$f_k(x_k)=\frac{1}{2}||A_kx_k-b_k||^2$ where $A_k\in\Rset^{p\times n}$ and $b_k\in\Rset^p$, $k=1,\ldots,K$. Then, the subproblems in each iteration can be written as
\be\label{subprob:exchange_problem}
x_k^{\nu+1}=\argmin_{x_k}\frac{1}{2}||A_kx_k-b_k||^2+r||x_k-d_k^\nu||^2,\quad\forall k=1,\ldots,K,
\ee
for some $r\in\Rset_+$ and $d_k^\nu\in\Rset^n$. Notice that the matrices $A_k^TA_k+2r\textbf{I}_n,k=1,\ldots,K$ are positive definite since $r>0$. We only have to compute $(A_k^TA_k+2r\textbf{I}_n)^{-1}$ for one time before the iterations start.  In the experiments, we randomly generate the optimal solution $x_k^*,k=1,\ldots,K-1$ by the standard normal distribution and set $x_K^*=-\sum_{k=1}^{K-1}x_k^*$. The matrices $A_k,k=1,\ldots,K$ are generated from standard  Gaussian distribution and we let $b_k=A_kx_k^*$. In this setting, $x^*$ is an optimal solution to \eqref{exchange problem} but not necessarily the unique one, and the optimal value is 0. We set  $K=20$, $n=1000$, $p=800$ , and none of  $f_k(x_k),k=1,\ldots,K$ is strongly convex.  We compare the performance of ADA with VSADMM and Prox-JADMM mentioned in Section \ref{Relation: ADMM}. For the implementation of VSADMM and Prox-JADMM, we use codes provided in \cite{deng2017parallel}.  For the proximal parameters of ADA, we set $(\rho,c)=(10,10)$ in the experiment. 

\begin{figure}[h]
	\centering
	\includegraphics[width=0.8\linewidth,height=0.25\textheight]{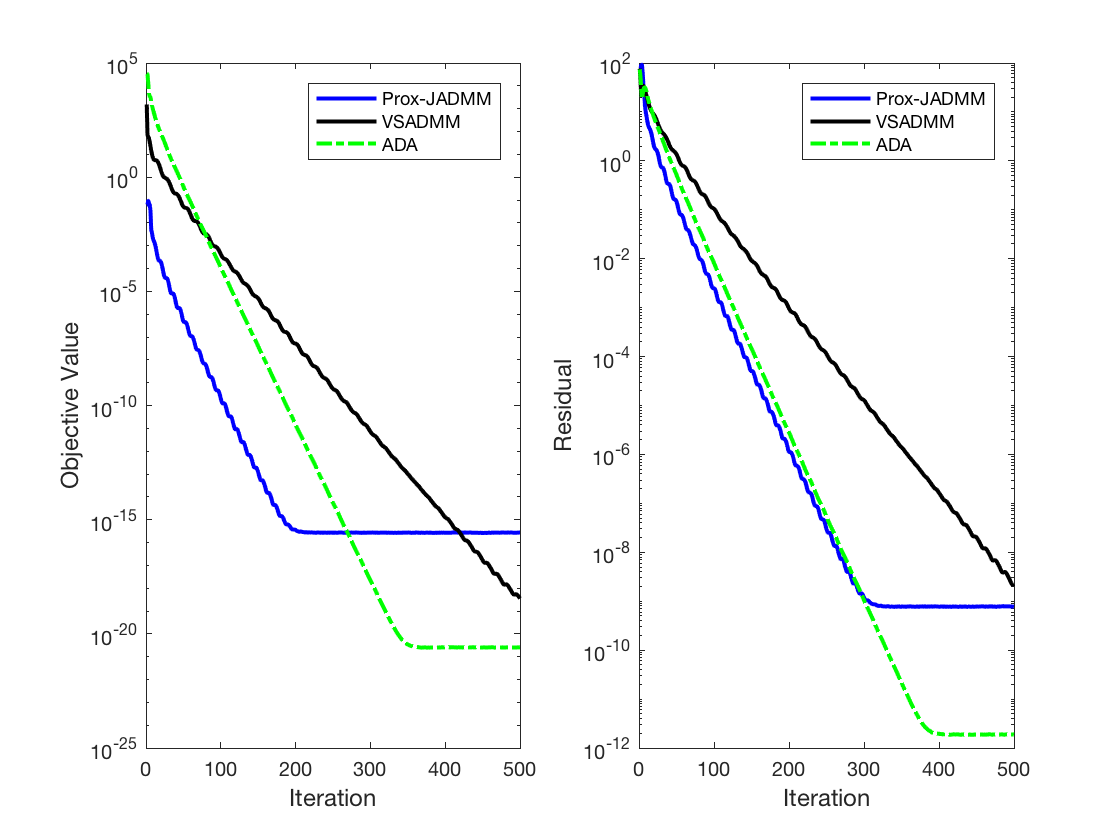}
	\caption{Exchange Problem: $K = 20, n = 1000, p = 800$. Convergence results versus iteration. }
	\label{fig:exchange_problem}
\end{figure}

\begin{figure}[h]
	\centering
	\includegraphics[width=0.8\linewidth,height=0.25\textheight]{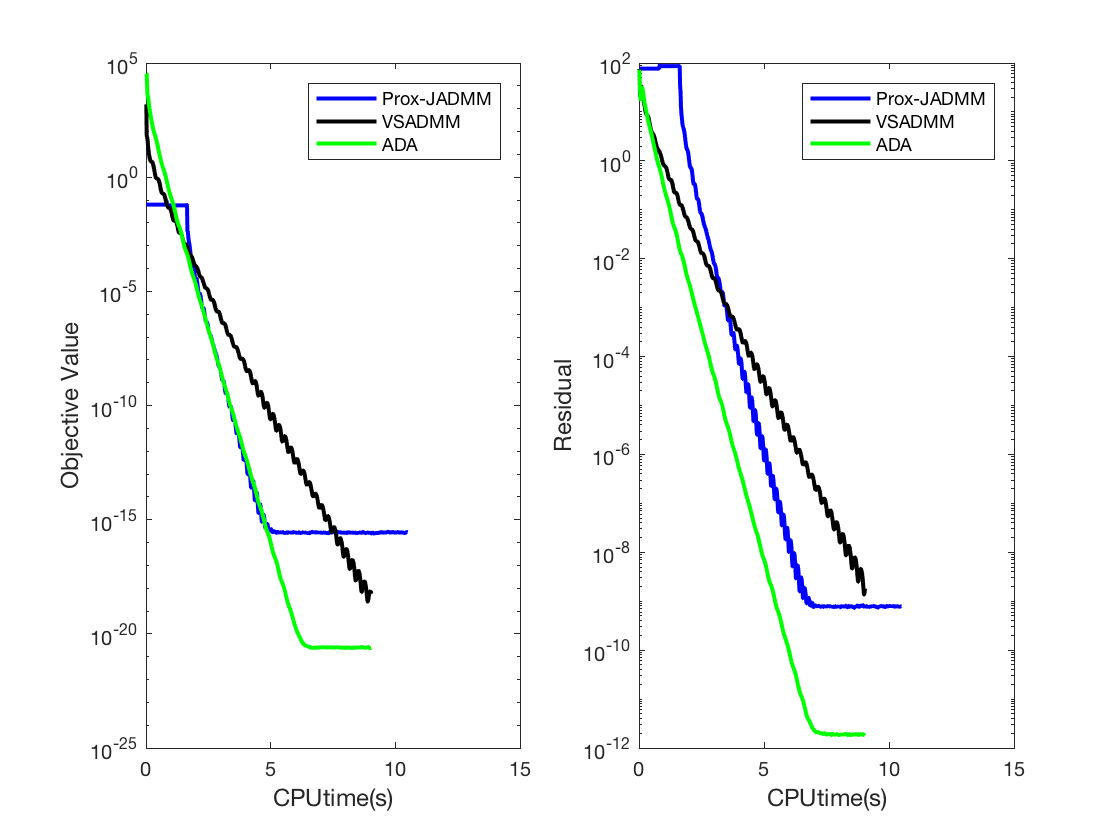}
	\caption{Exchange Problem: $K = 20, n = 1000, p = 800$. Convergence results versus time.}
	\label{fig:exchange_problem_time}
\end{figure}

For all of the  algorithms, we start from the same initial point $(x^0,y^0)=(\textbf{0},\textbf{0})$ and  run 500 iterations.  Figure \ref{fig:exchange_problem} shows the objective function value $\sum_{k=1}^{K}f_k(x_k)$ and the  residual $||\sum_{k=1}^{K}x_k||$ of each iteration for the average outcome of 10 random simulations. We can see that ADA shows a better convergence of the objective value compared with VSADMM and is slower than Prox-JADMM in terms of iterations. However, Prox-JADMM requires extra computational time to update the proximal parameters which is shown in Figure \ref{fig:exchange_problem_time}. Overall, ADA shows competitive convergence results in this experiment compared with two variants of the classical ADMM method which facilitate parallelization.

\subsection{Distributed $\ell_1$-regularized logistic regression}
Here, we use  iADA to solve the convex regularization problem \eqref{conv_reg_v} with a modest number of features but a relative large number of training examples. Many statistical problems belong to this regime, with a large $n$ and a small $d$ dataset. In particular, we consider the following $\ell_1$-regularized logistic regression:
\begin{equation}\label{l1-logistic}
\begin{aligned}
\min_{x\in\Rset^d} F(x)=\sum_{j=1}^{n}\ell(x;(a_j,b_j))+\lambda ||x||_1\\
\end{aligned}
\end{equation}
where $(a_j,b_j)\in\Rset^{d+1},j=1,\ldots,n$ and $\ell(x;(a_j,b_j))=\log(1+\exp(-b_ja_j^Tx))$. For the purpose of parallel computation, we partition $A\in\Rset^{n\times d}$ and $b\in\Rset^n$ into $N$ blocks
\[A=\begin{bmatrix}
A_1\\\vdots\\A_N
\end{bmatrix}\quad \text{and}\quad 
b=\begin{bmatrix}
b^1\\\vdots\\b^N
\end{bmatrix},
\]
with $A_i\in\Rset^{n_i\times d}$ and $b^i\in\Rset^{n_i}$. Define $\bar{n}_i=\sum_{j=1}^{i}n_j$  and we notice $\bar{n}_0=0$ and $\bar{n}_N=\sum_{j=1}^{N}n_j=n$. By introducing variables $x_i\in\Rset^d,i=1,\ldots,N$, \eqref{l1-logistic} can be transformed into the following:
\begin{equation}\label{distr:loss}
\begin{aligned}
 \min_{x_i,z\in\Rset^d}& \sum_{i=1}^{N}\ell_i(x_i;(A_i,b^i))+\lambda||z||_1\\
 \text{    s.t.      }&\quad  x_i-z=0, \quad i=1,\ldots,N.
\end{aligned}
\end{equation}
where $\ell_i(x_i;(A_i,b^i))=\sum_{j=\bar{n}_{i-1}+1}^{\bar{n}_i}\log(1+\exp(-b_ja_j^Tx_i))$.  In our experiment, we use two publicly available datasets: (1) the \texttt{w8a} dataset (49749 examples and 300 features) and (2) the \texttt{ijcnn1}  dataset (49990 examples and 22 feature). 
The main step of  iADA algorithm is given by  
\begin{equation}\label{distr_LR:iADA}
\begin{aligned}
&x_i^{\nu+1}\approx \argmin_{x_i}\underbrace{\ell_i(x_i;(A_i,b^i))+\frac{\rho}{4}||x_i-w_{x,i}^\nu+\frac{2}{\rho}y_{x,i}^\nu||_2^2+\frac{1}{2c}||x_i-x_i^\nu||_2^2}_{\phi_{i,\rho,c}^\nu(x_i)},\\
&z^{\nu+1} = \argmin_{z}\lambda_1||z||_1+\frac{\rho}{4}\sum_{i=1}^{N}||z+w_{z,i}^\nu-\frac{2}{\rho}y_{z,i}^\nu||_2^2+\frac{1}{2c}||z-z^\nu||_2^2,
\end{aligned}
\end{equation}
where $w_x^\nu=(w_{x,1}^\nu,\ldots,w_{x,N}^\nu)\in\Rset^{Nd},y_x^\nu=(y_{x,1}^\nu,\ldots,y_{x,N}^\nu)\in\Rset^{Nd}$ and $w_z^\nu,y_z^\nu\in\Rset^{Nd}$. The $x_i$ update involves an $\ell_2$ regularized logistic regression which cannot be solved exactly. Here, we use the L-BFGS algorithm to solve them until the inexact criteria \eqref{criteria_A} and \eqref{criteria_B} are satisfied. Such criteria can be checked by identifying the norm of the gradient $||\nabla\phi_{i,\rho,c}^\nu(x_i)||$. For the $z$ update, exact solutions can be derived by the soft threshold operator.

For comparison, we consider the inexact ADMM (iADMM) method proposed in \cite{eckstein1992douglas,he2002new}.  Similar subproblems as \eqref{distr_LR:iADA} will arise for $x_i$ and $z$ updates. An analogous inexact criterion as  \eqref{criteria_A} are proposed in \cite{eckstein1992douglas,he2002new} to guarantee the convergence of the inexact ADMM and can also be verified by examining the norm of the gradient in the $x_i$ updates.

In the experiment, we set $\epsilon_\nu=\frac{1}{\nu^\gamma}$ with $\gamma=1.0,1.5,2.0$ to control the inexactness of the $x_i$ updates in both algorithms. We also consider different partitions with $N=20, 50$. For the implementation of  iADMM, we use a  step-length 1.618 and a fixed penalty parameter 10  after tuning. For iADA, we choose the proximal parameters $(\rho,c)=(10,10)$. Both algorithms terminated when 
\[
\frac{\sum_{i=1}^{N}||x_i^\nu-z^\nu||_2}{N||z^\nu||_2}\le10^{-6}\text{ and }\frac{|F(z^{\nu})-F(z^*)|}{\max\{1,|F(z^*)|\}}\le 10^{-10}
\]
are satisfied. $F(z^*)$ is the optimal solution of \eqref{distr:loss} derived by running iADMM for 2000 iterations.

The computational results  are presented in Table \ref{table:distr-loss}. The datasets are listed in the first column. The numbers of partitions $N$ and the inexactness parameter $\gamma$ are given in columns two and three separately. The $\infty$ symbol in the third column represents the exact $x_i$ updates achieved by setting $\epsilon_\nu=1e-10$ in all iterations. The average number of iterations (upon round off) for iADA and iADMM are given in the next two columns. The total amount of L-BFGS updates for both methods are presented in columns 6-7 and the average CPU time (in seconds) for these methods are given in the last two columns.
\begin{table}[h]
	\caption{Comparison of iADA and iADMM for solving \eqref{distr:loss}.}\label{table:distr-loss}
	\centering
\begin{tabular}{|c|c|c|c c|c c|c c|}
	\hline
	\multirow{2}{*}{Dataset}& \multirow{2}{*}{N} &  \multirow{2}{*}{$\gamma$} &  \multicolumn{2}{c|}{Iteration}   & \multicolumn{2}{c|}{L-BFGS}  &  \multicolumn{2}{c|}{CPU time}  \\ 
	&  &  &  iADA&  iADMM& iADA &iADMM  & iADA  &  iADMM\\ 
	\hline 
	\multirow{8}{*}{w8a}	& \multirow{4}{*}{20} & 1.0 & 274  & 380 &  70361 & 83703 &  \textbf{24.00} & 30.00 \\ \cline{3-3}
	&  & 1.5 & 169 & 197 &41089 &   58199 &  \textbf{14.18} &  19.15\\ \cline{3-3}
	&  & 2.0 &  164 & 195 & 44616 &  65647  & \textbf{15.02}   & 20.87  \\ \cline{3-3}
	&  & $\infty$ & 150 & 133 &49945 & 54928 & \textbf{15.70} &17.78 \\\cline{2-9}
	& \multirow{4}{*}{50} & 1.0 & 172 &  211 & 88460 & 127538  & \textbf{16.00} & 22.34 \\ \cline{3-3}
	&  &  1.5& 140 & 120 &78594  & 72673  & 13.10& \textbf{10.66} \\ \cline{3-3}
	&  & 2.0 & 99 & 88 & 67909 &  60368 & 11.75 &  \textbf{10.10} \\\cline{3-3}
	&  & $\infty$ & 106 & 70 &77627 & 62893 &13.19 & \textbf{10.50} \\\hline
		\multirow{8}{*}{ijcnn1}& \multirow{4}{*}{20} &  1.0 & 202 &  276& 49378 &79120  &  \textbf{17.02} & 26.80 \\ \cline{3-3}
	&  & 1.5 & 114 & 135 & 29741 & 42142 &  \textbf{10.16} & 14.10 \\ \cline{3-3}	
	&  & 2.0 & 112 &134  & 31308 & 46742 &  \textbf{10.50} & 15.11 \\ \cline{3-3}
	&  & $\infty$ & 190 & 186 &73111 & 73193 & 23.02 & \textbf{22.15} \\\cline{2-9}
& \multirow{4}{*}{50} &  1.0 & 106 &228  & 68001 & 115891 & \textbf{11.74} & 22.05\\ \cline{3-3}
&  & 1.5 & 107 & 112 & 58093 &64195  &\textbf{10.97}  & 11.58 \\ \cline{3-3}
&  & 2.0 & 99 & 88 & 57777 & 50099 & 10.69 & \textbf{9.27}  \\ \cline{3-3}
&  & $\infty$ & 95 & 83 &68652 & 69291 & 11.32 & \textbf{11.21} \\
	\hline 
\end{tabular} 
\end{table}

From Table \ref{table:distr-loss}, we see that when $\gamma=1.5\text{ or }2.0$, iADA shows better performance in the case $N=20$ while iADMM converges faster when $N=50$. For both algorithms, the CPU time is much longer in the case of $\gamma=1.0$ when the convergence is not guaranteed in theory. Finally, compared with the exact update, it takes more iterations for the inexact version of both algorithms to converge but with shorter CPU time. This phenomenon results from the large number of  L-BFGS updates in each iteration of exact ADA and ADMM.

\section{Conclusions}
In this paper, we study the convergence results of the ADA and its inexact version, the iADA, for solving multi-block separable convex minimization problems subject to linear constraints. First, we prove the global convergence and the $o(1/\nu)$ rate for the exact ADA when there exists a saddle point for the corresponding Lagrangian function.
Next, global convergence and local linear convergence for the iADA are  established under some mild assumptions and certain approximation criteria. 

Before ending this paper, we would like to discuss two possible directions related to the ADA. Firstly, we notice that both the primal PPA \cite{guler1992new} and the Augmented Lagrangian Method \cite{he2010acceleration} can be accelerated by utilizing the idea from Nesterov's seminal work \cite{nesterov1983method}. It is natural to ask whether we can accelerate the ADA based on similar techniques since all of them belong to the general PPA framework. Secondly, the applicability of the approximation criteria in \eqref{criteria_A} and \eqref{criteria_B} is limited in practice due to the summable requirement and  more implementable approximation criteria are needed  for practical problems.

\section*{Acknowledgments}
The authors are grateful to Professor R. Tyrrell Rockafellar for suggestions on this research project. Shu Lu's research is  supported by National Science Foundation under the grant DMS-1407241.

\bibliographystyle{spmpsci}
\bibliography{my_ref}

\end{document}